\documentclass[11pt,twoside]{article}
\usepackage{latexsym}
\usepackage{amssymb,amsbsy,amsmath,amsfonts,amssymb,amscd}
\usepackage{dsfont}
\usepackage{graphicx,color}
\usepackage[colorlinks=true]{hyperref}
\usepackage{epsfig,color}
\usepackage{dsfont}
\usepackage{algorithm}

\usepackage[amsmath,framed,thmmarks]{ntheorem}

\usepackage{placeins}
\usepackage{subcaption}
\usepackage{amsmath} 
\numberwithin{equation}{section}
\DeclareMathOperator{\sign}{sign}
\usepackage{longtable, array, booktabs}
 
\usepackage{tabularx}
\usepackage{mathtools}
\newcommand*\diff{\mathop{}\!\mathrm{d}}

\theoremstyle{plain}
\newtheorem{theorem}{Theorem}[section]

\newtheorem{proposition}[theorem]{Proposition}

\theoremstyle{nonumberplain}
\newtheorem{proof}{Proof.-}

\newcommand{\brho}{\boldsymbol{\rho}}
\definecolor{myblue}{rgb}{0.0352,0.4981,0.6509}

\newcommand{\e}{\epsilon}

\title{An asymptotic preserving scheme for kinetic equations describing bacterial traveling pulses}
\author{Canan Akkoyunlu\thanks{Department of Mathematics and Computer Science, İstanbul Kültür
University, 34158 Bakırköy, İstanbul, Türkiye (c.kaya@iku.edu.tr)} \and Jose A. Carrillo\thanks{Mathematical Institute, University of Oxford, Woodstock Road, Oxford, OX2 6GG, UK (jose.carrillo@maths.ox.ac.uk)}\and Gissell Estrada-Rodriguez\thanks{Universitat Politecnica de Catalunya, Barcelona, 08034, Spain. Centre de Recerca Matematica, 08193 Bellaterra, Barcelona, Spain (gissell.estrada@upc.edu)} \and  Xinran Ruan\thanks{School of Mathematical Sciences, Capital Normal University, Beijing, China, 100048 (xinran.ruan@cnu.edu.cn). 
} 
}
\date{}

\begin{document}

\maketitle
\begin{abstract}
We consider a kinetic model for bacterial chemotaxis that describes the propagation of traveling pulses driven by a nutrient and a self-produced chemoattractant. In the diffusive regime, the model converges to a macroscopic drift–diffusion system with bounded chemotactic fluxes. We construct a fully discrete asymptotic preserving (AP) finite difference scheme based on a micro–macro decomposition of the kinetic equation. The method is stable uniformly with respect to the scaling parameter and, as this parameter tends to zero, reduces to a consistent discretization of the limiting macroscopic equations. We prove positivity and mass conservation for the limiting scheme, and use a positivity-preserving correction for the kinetic solver. Numerical simulations reproduce the experimentally observed traveling pulses at both the kinetic and the macroscopic level, validate the theoretical traveling wave speed, and illustrate the convergence of the kinetic model to its macroscopic limit.
\end{abstract}

\section{Introduction}

Chemotaxis is the directed movement of organisms in response to chemical fluctuations in the environment. This behavior is essential for survival \cite{franca2006navigating} and plays a crucial role in many biological processes in multicellular organisms, such as wound healing, fertilization and inflammation \cite{rickert2000leukocytes},  while unicellular organisms use it as a foraging mechanism \cite{franca2010g,willard2006signaling}. Many bacteria, in particular, navigate toward favorable conditions by sensing gradients of nutrients or chemoattractants, and away from harmful substances.

Although bacterial chemotaxis has been studied since the end of the 19th century, the first precise mathematical model was proposed by Keller and Segel \cite{keller1971model}. They introduced a coupled system of partial differential equations (PDEs) for the cell density $\rho(t,x)$ and the chemical concentration $S(t,x)$:
\begin{subequations}\label{eq:chemicals}
\begin{align}
\frac{\partial \rho}{\partial t} &= \nabla\cdot (D_\rho (\rho, S) \nabla \rho - \chi (\rho, S) \nabla S) + f(\rho, S),
  \label{eq: KS1}\\
\frac{\partial S}{\partial t} &= D_S \Delta S + g (\rho, S) - h(\rho, S)S.
  \label{eq: KS2}
\end{align}
\end{subequations}
Here $D_\rho(\rho,S)$ denotes the bacterial diffusivity, $\chi(\rho,S)$ the chemotactic
sensitivity, and $f$, $g$, and $h$ describe cell growth and death, chemical production,
and degradation, respectively.

We consider the migration of \textit{E.~coli} in response to a nutrient $N$ and a self-produced chemoattractant $S$, motivated by the microchannel experiments in \cite{saragosti2010mathematical}. In these experiments, approximately $10^5$ cells are initially confined near the left end of a $2\,\mathrm{cm}$ channel, while nutrient is supplied from the right. Within about $100\,\mathrm{s}$, the bacteria organize into a dense pulse that propagates toward the nutrient source with nearly constant speed and mass. The pulse exhibits a characteristic asymmetric profile, with a diffuse front and a steep rear, which is maintained during propagation. Since the duration of the experiment is short compared with the cell division cycle, bacterial growth is neglected.

Following the kinetic descriptions of Alt \cite{alt1980biased} and Dunbar and Othmer
\cite{othmer1988models}, we model the bacterial population by the mesoscopic density
$f(t,x,v)$, where $x \in \mathbb{R}$ and $v \in V$ denotes the velocity. Since we are
interested in planar traveling waves, the dynamics are governed by the one-dimensional
kinetic equation
\begin{equation}
    \partial_t f(t,x,v)+{v}\partial_x f(t,x,v)= \int_{V}\Bigl(T(t,x,v,v')f(t,x,v')-  T(t,x,{v}',v)f(t,x,v)\Bigr)\diff v'\ .\label{eq: main kinetic}
\end{equation}
 The turning kernel $T$ describes the change in direction during the tumbling phase of the movement. $T$ depends on the space and time variables through the chemical signals in the medium,  and $T(\cdot, v,v')$ represents the change of velocity from $v'$ to $v$. For simplicity we are going to assume that $T$ only depends on the current velocity with the condition that $\int_VT(t,x,v)\diff v=1$.
 
 Following \cite{saragosti2010mathematical}, we consider the influence of two chemical species: a chemotactic signal $S(t,x)$ released by the bacteria, and a nutriet $N(t,x)$ consumed by the them. We assume that both chemical cues contribute additively to the turning frequency. 
More precisely, we define
\[
\widetilde{\psi}_1(t,x,v)
=
\chi_S \phi\!\left(\frac{DS}{Dt}\Big|_v\right)
+
\chi_N \phi\!\left(\frac{DN}{Dt}\Big|_v\right),
\qquad
\frac{D}{Dt} = \partial_t + v \partial_x,
\]
where $\phi$ is the intracellular response function and $\chi_S$, $\chi_N$ are the
chemotactic sensitivities. Since the velocity average of $\widetilde{\psi}_1$ does not
vanish in general, we introduce the centered perturbation
\begin{equation}\label{eq: centered psi1}
\psi_1(t,x,v)
=
\widetilde{\psi}_1(t,x,v)
-
\frac1{|V|}\int_V \widetilde{\psi}_1(t,x,w)\,dw,
\end{equation}
which satisfies $\int_V \psi_1(t,x,v)\,dv = 0.$
The rescaled turning kernel is then defined by
\[
T_\varepsilon(t,x,v)
=
|V|\,
\frac{1+\varepsilon \widetilde{\psi}_1(t,x,v)}
{\displaystyle \int_V \bigl(1+\varepsilon \widetilde{\psi}_1(t,x,w)\bigr)\,dw},
\]
and admits the first-order expansion
\begin{equation}\label{eq: T expansion}
T_\varepsilon(t,x,v)=1+\varepsilon \psi_1(t,x,v)+O(\varepsilon^2),
\end{equation}
with
$\int_V T_\varepsilon(t,x,v)\,dv = |V|.$

The chemoattractant and nutrient concentrations satisfy the reaction--diffusion equations
\begin{subequations}
\label{eq:SN_intro}
\begin{align}
\partial_t S &= D_S \partial_{xx} S - \alpha S + \beta \rho, \\
\partial_t N &= D_N \partial_{xx} N - \gamma \rho N,
\end{align}
\end{subequations}
where
\[
\rho(t,x) = \int_V f(t,x,v)\,dv\ ,
\]
denotes the macroscopic bacterial density.
Here the chemoattractant is produced by the bacteria at rate $\beta$ and naturally degraded at rate $\alpha$, while the nutrient is consumed at rate $\gamma$.

The main objective of this paper is to construct and validate a numerical scheme that accurately captures the dynamics of the kinetic model for all values of the scaling parameter $\varepsilon$, and recovers the corresponding macroscopic drift–diffusion system in the diffusive limit $\varepsilon\to0$. In particular, we investigate whether the traveling pulses observed at the macroscopic level are consistently reproduced by the kinetic description. To this end, we develop an asymptotic preserving (AP) scheme \cite{filbet2010class,jin1999efficient}, which remains stable as $\varepsilon\to0$ and, for fixed discretization parameters, automatically reduces to a consistent discretization of the limiting macroscopic equations. Our approach builds on the micro–macro decomposition introduced in \cite{lemou2008new} (see Section~\ref{sec: micro macro decomp}). Other AP approaches for kinetic chemotaxis models can be found in \cite{CY13,LWZ18,CKLO19} and the references therein.

In the diffusive limit, the kinetic model formally converges to the macroscopic system
\begin{subequations}
\label{basicmodel}
\begin{align}
\partial_t \rho
&= D_\rho \Delta \rho
   - \nabla \cdot \bigl( \rho u_S + \rho u_N \bigr), \\
\partial_t S
&= D_S \Delta S - \alpha S + \beta \rho, \\
\partial_t N
&= D_N \Delta N - \gamma \rho N,
\end{align}
\end{subequations}
where the drift velocities $u_S$ and $u_N$ are derived from the microscopic response function and take the form
\begin{equation}
\label{eq:flux_intro}
u_S
=
\chi_S J_\phi(-\varepsilon \partial_t S, |\nabla S|)
\frac{\nabla S}{|\nabla S|},
\qquad
u_N
=
\chi_N J_\phi(-\varepsilon \partial_t N, |\nabla N|)
\frac{\nabla N}{|\nabla N|}.
\end{equation}
 Here $\phi$ is the so-called response function
 (Section~\ref{sec: maro system}), which
describes how a single bacterium reacts to changes in the concentration of the
chemoattractant $S$ (or of the nutrient $N$) along its trajectory. The function $J_\phi$
encodes, at the macroscopic level, these microscopic features of the individual response.
The parameter $\varepsilon$ is the ratio between the pulse speed and the speed of individual
bacteria; in the experiments $\varepsilon \approx 0.1$, so the pulse travels one order of
magnitude slower than each individual.

In contrast with the classical Keller–Segel fluxes, these drift terms remain uniformly bounded. This provides a more realistic description of chemotactic motion and allows an analytical characterization of traveling pulses \cite{calvez2019chemotactic}: the flux $u_S$ aggregates the cells into a pulse, while $u_N$ drives this pulse toward higher nutrient concentrations. In the stiff-response limit, where $\phi$ is a step function, the chemotactic drift reduces to $\chi\,\mathrm{sign}(\partial_x S)$. For this case, in \cite{calvez2023nonlinear} the authors studied the one-dimensional model without nutrient and with $S$ in quasi-stationary equilibrium, and proved the nonlinear stability of its stationary states, which are exponentially decaying clusters, i.e. zero-speed pulses.

The paper is organized as follows. Section~\ref{sec: maro system} derives the macroscopic drift–diffusion system from the underlying kinetic chemotaxis model and recalls the theoretical traveling wave speed. Section~\ref{sec: micro macro decomp} introduces the micro–macro decomposition used to construct the asymptotic preserving formulation. In Section~\ref{sec:scheme}, a fully discrete AP finite difference scheme is developed, and its asymptotic preserving, positivity, and mass conservation properties are analyzed. Section~\ref{sec: numerics} presents numerical experiments validating the traveling pulse dynamics, the dependence of the wave speed on the chemotactic sensitivities, and the convergence of the kinetic model towards its macroscopic limit as the scaling parameter tends to zero. Finally, Section~\ref{sec: conclusions} summarizes the main results and discusses possible directions for future work.

\section{Macroscopic system of equations}\label{sec: maro system}
In this section we derive the macroscopic system that is obtained from \eqref{eq: main kinetic} following \cite{saragosti2010mathematical}. We introduce two relevant parameters: the maximum speed of a single bacterium $c = \max\{|v|: v \in V\}$ and the mean turning frequency $\lambda_0$.
The main unknown is the speed of the traveling pulse, denoted by $\sigma$. We start by rescaling the kinetic model (\ref{eq: main kinetic}) into a non-dimensional form using the non-dimensional variables, denoted with tildes, as
\begin{equation}
    \tilde{t}= \frac{{t}}{\bar{t}}\ , \quad \tilde{x} = \frac{x}{\bar{x}}\ , \quad {\tilde{V} = \frac{V}{c}}\ , \quad \tilde{T}=\frac{T}{\lambda_0}\ ,
    \nonumber
\end{equation}
to model the traveling wave with $\bar{x} = \sigma\bar{t}$. According to the experiments, the pulse
speed $\sigma$ is much lower than the speed of a single bacterium $c$, and this difference is measured by the ratio $\varepsilon = \sigma / c$, where $\varepsilon \ll 1$.

The kinetic equation can be written in the following form, where we  drop the tilde for simplicity of notation 
{
\begin{equation}
    \varepsilon \partial_t f^\varepsilon+{v}\partial_x f^\varepsilon = \frac{\mu}{\varepsilon}\Bigl[\int_V T_\varepsilon(t,x,v')f^\varepsilon(t,x,v')\diff v'-|V|T_\varepsilon(t,x,v)f^\varepsilon(t,x,v)\Bigr]\ , \label{eq: inital model with scaling}
\end{equation}}where $\mu/\varepsilon=\lambda_0 \bar{x}/c$ and with velocity space \( V = [-1,1] \) (one-dimensional case). 
Here 
\begin{equation}\label{eq: T app}
T_\varepsilon(t,x,v)
=
1+\varepsilon\psi_1(t,x,v).
\end{equation} 
is the first-order approximation of \eqref{eq: T expansion}
{and \(\psi_1\) denotes the centered first-order perturbation defined in \eqref{eq: centered psi1}.
We consider a regular expansion of $f^\varepsilon$ in $\varepsilon$ up to second order
\begin{equation*}
f^\e(t,{x},{v})=f^0(t,{x},{v})+\varepsilon f^1(t,{x},{v})+\varepsilon^2 f^2(t,{x},{v})+\mathcal{O}(\varepsilon^3)\ , 
\end{equation*}
where $\int_V f^\e(t,{x},{v})\diff {v} = \int_V f^0(t,{x},{v})\diff {v} = \rho(t,{x})$, therefore $\int_V f^i(t,{x},{v})\diff {v}=0$, $\forall i\geq 1$.

Substituting this ansatz into \eqref{eq: inital model with scaling} and collecting powers of $\varepsilon$
\iffalse
{\begin{align*}
    \varepsilon^2\partial_t f^0&+\varepsilon{v}\partial_x\big(f^0+\varepsilon f^1  \big) 
    =-\psi_0(v) f^0+\int_V\psi_0(v')f^0(v')\diff v'\\
    &+\varepsilon\int_V\Bigl(\psi_0(v')f^1(v')+\psi_1(t,x,v')f^0(v')-\psi_0(v)f^1(v)-\psi_1(t,x,v)f^0(v) \Bigr)\diff v'+{\mathcal{O}(\varepsilon^3)}
      \ .
\end{align*}}

\noindent Identifying the different equations in powers of $\varepsilon$, we obtain
\fi
\begin{align}
    \varepsilon^0: & \quad f^0(t,{x},{v})=\frac{\rho(t,{x})}{|V|}\ ,\label{eq: varepsilon0}\\
    \varepsilon^1: & \quad 
    v \partial_x f^0
=
\mu \left(
\int_V \psi_1(v') f^0(v')\,dv'
- |V| \psi_1(v) f^0(v)
- |V| f^1(v)
\right)\ ,\label{eq: varepsilon1}\\
    \varepsilon^2: & \quad \partial_tf^0+{v}\partial_x f^1 =\mu\Bigl(\int_V\psi_1(v')f^1(v')\diff v'-|V|\psi_1(v)f^1(v)-|V|f^2(v)\Bigr) \ \label{eq: varepsilon2}.\
\end{align}
Next, after replacing $f^0$ given by  \eqref{eq: varepsilon0}, we multiply \eqref{eq: varepsilon1} by ${v}$ and we integrate again with respect to ${v}$ to obtain
\begin{align}
    \mu|V|\int_V vf^1(t,x,v) \diff v=&-\partial_x\frac{\rho(t,x)}{|V|}\int_V v^2\diff v-\mu\rho(t,x)\chi_S\int_V\phi(\varepsilon\partial_tS+v\partial_xS)v\diff v\nonumber\\ &-\mu\rho(t,x)\chi_N\int_V\phi(\varepsilon\partial_tN+v\partial_xN)v\diff v.\label{eq: mean direction}
\end{align}
We define
\begin{equation}
\label{eq: moment}
\begin{gathered}
D_\rho
=
\frac{1}{\mu |V|^2}\int_V v^2\,\diff v,
\qquad
u_S
=
-\frac{\chi_S}{|V|}
\int_V \phi(\varepsilon\partial_t S+v\partial_x S)\,v\,\diff v,
\\[4pt]
u_N
=
-\frac{\chi_N}{|V|}
\int_V \phi(\varepsilon\partial_t N+v\partial_x N)\,v\,\diff v .
\end{gathered}
\end{equation}
In the limiting case where the internal response function is $\phi(Y)=-\phi_0\sign(Y),$
the flux terms reduce, up to the constant factor determined by
$\phi_0$ and the velocity measure, to \cite{calvez2019chemotactic}
\[
u_S=\chi_S\sign(\partial_x S),
\qquad
u_N=\chi_N\sign(\partial_x N).
\]
Finally, integrating \eqref{eq: varepsilon2} in ${v}$ and using \eqref{eq: moment} we obtain
\[
\partial_t\rho=D_\rho\partial_{xx}\rho-\partial_x(\rho u_S+\rho u_N)\ .
\]
{
The diffusion coefficient is determined by the velocity distribution and
the turning frequency as in \eqref{eq: moment}. 
In particular, when the normalized scaling \(\mu |V|=1\) is used, this
coefficient becomes
\begin{equation}
\label{def:Drho}
  D_\rho
=
\frac1{|V|}
\int_V v^2\,dv .  
\end{equation}
For example, when \( V = [-1,1] \), we have $D_\rho = \frac{1}{3}$.
}
% In particular, since \( V = [-1,1] \), we have $|V|=2$ and $
% \int_{-1}^1 v^2 dv = \frac{2}{3}.
% $
% Hence, $D_\rho=
% \frac{1}{6\mu}$ and if, for instance, \( \mu = 1 \),
% $D_\rho = \frac{1}{6}$.
We refer to \cite{emako2015existence,emako2016traveling} for more details of the derivation.
Together with \eqref{eq:SN_intro} for $S$ and $N$, this system gives a traveling wave solution where the speed is given by the following theorem.

\begin{theorem}[\cite{calvez2019chemotactic}]\label{thmacro}
There exist a speed $\sigma>0$, and two positive values $N_- < N_+$, such that the system \eqref{basicmodel} admits traveling wave solutions ($\rho (z), S(z), N(z)$) for $z=x - \sigma t$ with the following boundary conditions:
\begin{equation}
    \lim_{|z| \to \infty} \rho = \lim_{|z| \to \infty} S = 0, \hspace{0.6cm} \lim_{z \to -\infty} N = N_- < N_+ = \lim_{z \to + \infty} N.
    \nonumber
\end{equation}
Assume in addition that $\partial_z S$ changes sign only once, and that $\partial_z N$ is positive everywhere. Then, the speed $\sigma > 0$ is uniquely determined by the following implicit relation,
\[
\chi_N - \sigma = \chi_S \frac{\sigma}{\sqrt{4 D_S \alpha + \sigma^2}}.
\]
\end{theorem}
The proof in \cite{calvez2019chemotactic} relies on the ansatz that the chemoattractant $S$ has a unique maximum at $z=0$ and that the nutrient $N$ is increasing. This ansatz decouples system \eqref{basicmodel}, so that one can solve explicitly for $\rho$ and then for $S$; validating the ansatz a posteriori yields the implicit relation for $\sigma$ stated in Theorem \ref{thmacro}.

{
\section{Micro-macro decomposition}
\label{sec: micro macro decomp}

Throughout this section, we use the scaling $\mu |V|=1$. 
To construct an asymptotic preserving scheme consistent with the diffusive limit, we introduce the (constant) equilibrium distribution 
\[
M(v)=\frac1{|V|},
\qquad \int_V M(v)\,dv=1.
\]
Since the leading-order term is given by \(f^0=\rho/|V|\), we decompose the kinetic density as
\begin{equation}\label{f_expansion}
f^\varepsilon(t,x,v)=\frac{\rho^\varepsilon(t,x)}{|V|}+\varepsilon g^\varepsilon(t,x,v),
\qquad
\int_V g^\varepsilon(t,x,v)\,dv=0.
\end{equation}

We denote by \(\langle \cdot \rangle\) the integration with respect to \(v\) over \(V\), namely
\[
\langle h\rangle = \int_V h(v)\,dv,
\]
and define the projection operator
$
\Pi h = \langle h\rangle M.
$ Then \((I-\Pi)g^\varepsilon=g^\varepsilon\), where $I$ is the identity.

Integrating \eqref{eq: inital model with scaling} with respect to \(v\), and using the conservation property of the $T_\varepsilon$ operator together with
$
\langle vM\rangle = \frac1{|V|}\int_V v\,dv = 0,
$
we obtain the macroscopic equation
\[
\partial_t \rho^\varepsilon + \partial_x \langle v g^\varepsilon\rangle = 0.
\]

% Next, we define 
% \[
% \mathcal{Q}_\varepsilon(f)
% =
% \int_V T_\varepsilon(t,x,v')f(v')\,dv'
% -
% |V|T_\varepsilon(t,x,v)f(v).
% \]
% % The normalized turning kernel admits the first-order expansion
% % \[
% % T_\varepsilon(t,x,v)
% % =
% % 1+\varepsilon\psi_1(t,x,v)+O(\varepsilon^2),
% % \qquad
% % \langle \psi_1\rangle=0.
% % \]
% In the following construction, we use the first-order approximation of \eqref{eq: T expansion} given by
% \[
% T_\varepsilon^{\rm app}(t,x,v)
% =
% 1+\varepsilon\psi_1(t,x,v).
% \]
% It preserves the normalization
% $
% \langle T_\varepsilon^{\rm app}\rangle=|V|.
% $
We define the operator
\[
\mathcal{Q}_\varepsilon(f)
=
\int_V T_\varepsilon(t,x,v')f(v')\,dv'
-
|V|T_\varepsilon(t,x,v)f(v),
\]
and using \eqref{f_expansion} we obtain
\[
(I-\Pi)\mathcal Q_\varepsilon(f^\varepsilon)
=
-\varepsilon\bigl(|V|g^\varepsilon+\rho^\varepsilon\psi_1\bigr)
-\varepsilon^2|V|(I-\Pi)(\psi_1g^\varepsilon).
\]
Applying \(I-\Pi\) {to the kinetic equation \eqref{eq: inital model with scaling}} with
\(\mathcal Q_\varepsilon\) gives
\[
\varepsilon^2 \partial_t g^\varepsilon
+
\varepsilon (I-\Pi)(v\partial_x g^\varepsilon)
+
vM\,\partial_x \rho^\varepsilon
=
-\mu\bigl(|V|g^\varepsilon+\rho^\varepsilon\psi_1\bigr)
-\mu\varepsilon |V|(I-\Pi)(\psi_1g^\varepsilon).
\]
Using \(M=1/|V|\) together with the normalization \(\mu |V|=1\), we obtain
the micro-macro system
\begin{equation}\label{eq:micro-macro-normalized}
\left\{
\begin{aligned}
&\partial_t \rho^\varepsilon
+\partial_x \langle v g^\varepsilon\rangle = 0,\\
&\partial_t g^\varepsilon
+\frac1\varepsilon (I-\Pi)(v\partial_x g^\varepsilon)
=
\frac1{\varepsilon^2}
\left(
-\frac{v}{|V|}\partial_x \rho^\varepsilon
-g^\varepsilon
-\frac{\rho^\varepsilon}{|V|}\psi_1
\right)
-\frac1\varepsilon (I-\Pi)(\psi_1g^\varepsilon).
\end{aligned}
\right.
\end{equation}

}

{
\section{An asymptotic preserving finite difference scheme}
\label{sec:scheme}

In this section we construct an asymptotic preserving finite difference
scheme for the micro-macro formulation \eqref{eq:micro-macro-normalized}. We present the method in one space dimension. The scheme is
designed to remain stable as \(\varepsilon\to0\) and to reduce to a
consistent discretization of the macroscopic drift-diffusion limit.

Let $\Omega = [x_{\min}, x_{\max}]$ be the computational domain and \(x_{j+\frac12}=x_{\min}+(j+\frac12)\Delta x\),
\(x_j=x_{\min}+j\Delta x\), and \(t_n=n\Delta t\). We approximate
\(\rho\), \(S\), and \(N\) at the half grid points, while the microscopic
unknown \(g\) is approximated at the integer grid points,
\[
\rho_{j+\frac12}^n\approx \rho(t_n,x_{j+\frac12}),\qquad
S_{j+\frac12}^n\approx S(t_n,x_{j+\frac12}),\qquad
N_{j+\frac12}^n\approx N(t_n,x_{j+\frac12}),
\]
and $
g_{j,k}^n\approx g(t_n,x_j,v_k),\ v_k\in V_h $ where \(V_h=\{v_k\}\) is a symmetric discrete velocity set satisfying 
$v_{-k}=-v_k$.

For a half-grid scalar variable \(q_{j+\frac12}^n\), we use
\[
\delta_t^+ q_{j+\frac12}^n
=
\frac{q_{j+\frac12}^{n+1}-q_{j+\frac12}^n}{\Delta t},
\ \ \
(\delta_x q)_j^n
=
\frac{q_{j+\frac12}^n-q_{j-\frac12}^n}{\Delta x},\ \ \ \textnormal{and}\ \ \
(\delta_x^2 q)_{j+\frac12}^n
=
\frac{
q_{j+\frac32}^n
-2q_{j+\frac12}^n
+q_{j-\frac12}^n
}{\Delta x^2}.
\]
Moreover, we define the discrete divergence for \(g_{j,k}^n\) by
\[
(\delta_x g)_{j+\frac12,k}^n
=
\frac{g_{j+1,k}^n-g_{j,k}^n}{\Delta x}.
\]
Whenever a half-grid variable is needed at an integer grid point, we use
the local average
\[
q_j^n=\frac12\left(q_{j-\frac12}^n+q_{j+\frac12}^n\right).
\]
The value \((\delta_x^2 q)_j^n\) is defined in the same way by averaging
the two neighboring half-grid values.

We define
\[
\langle \eta\rangle_h=\sum_{k}\omega_k\eta_k,
\qquad \omega_k>0,
\qquad |V|_h=\sum_{k}\omega_k,
\qquad \Pi_h\eta=\frac{\langle\eta\rangle_h}{|V|_h}.
\]
We assume $\langle v\rangle_h=0$. Then $I-\Pi_h$ is the
projection onto the discrete zero-average space.

\paragraph{Discrete response functions}

At each time step, the response functions are evaluated explicitly. 
For $q\in\{S,N\}$, we use the scaled derivative
\begin{equation}
\mathcal D_\varepsilon^v q
=\varepsilon\partial_tq+v\partial_xq.
\end{equation}
We approximate $\partial_tq$ using the equations for \(S\) and
\(N\) at time $t^n$. 
We define
\[
\mathcal L_{S,j}^n
=
 D_S(\delta_x^2 S)_j^n
-\alpha S_j^n
+\beta \rho_j^n,
\quad 
\mathcal L_{N,j}^n
=
 D_N(\delta_x^2 N)_j^n
-\gamma \rho_j^n N_j^n .
\]
The discrete response arguments are
\[
\mathcal D_{\varepsilon,h}^{v_k}S_j^n
=\varepsilon\mathcal L_{S,j}^n+v_k(\delta_xS)_j^n,
\qquad
\mathcal D_{\varepsilon,h}^{v_k}N_j^n
=\varepsilon\mathcal L_{N,j}^n+v_k(\delta_xN)_j^n.
\]
All quantities on the right-hand sides are known at time $t^n$.

The velocity average of \(\phi(a+bv)\) does not vanish in general, even
when \(\phi\) is odd. We therefore use the normalized discrete responses
\[
\psi_{S,j,k}^n
=\chi_S\left[
\phi\!\left(\mathcal D_{\varepsilon,h}^{v_k}S_j^n\right)
-\Pi_h\phi\!\left(\mathcal D_{\varepsilon,h}^{v_k}S_j^n\right)
\right],
\psi_{N,j,k}^n
=\chi_N\left[
\phi\!\left(\mathcal D_{\varepsilon,h}^{v_k}N_j^n\right)
-\Pi_h\phi\!\left(\mathcal D_{\varepsilon,h}^{v_k}N_j^n\right)
\right].
\]
We set
$\psi_{j,k}^n
=
\psi_{S,j,k}^n+\psi_{N,j,k}^n ,$
and by construction,
\[
\langle \psi_{S,j,k}^n\rangle_h=0,
\qquad
\langle \psi_{N,j,k}^n\rangle_h=0,
\qquad
\langle \psi_{j,k}^n\rangle_h=0.
\]
Consequently, the discrete turning kernel $
T_{\varepsilon,j,k}^n
=
1+\varepsilon \psi_{j,k}^n$
satisfies
$
\langle T_{\varepsilon,j,k}^n\rangle_h=|V|_h.$
We also introduce
$
B_{S,j}^n
=
\langle v_k\psi_{S,j,k}^n\rangle_h,
\
B_{N,j}^n
=
\langle v_k\psi_{N,j,k}^n\rangle_h.
$
Then, the corresponding discrete drift velocities are defined by
\[
u_{S,j}^n
=
-\frac{B_{S,j}^n}{|V|_h},
\qquad
u_{N,j}^n
=
-\frac{B_{N,j}^n}{|V|_h}.
\]

\paragraph{Fully discrete micro-macro scheme}
The fully discrete AP scheme is given by
\begin{subequations}\label{eq:modified-main-scheme}
\begin{align}
&
\delta_t^+\rho_{j+\frac12}^n
+
\left\langle
v_k(\delta_x g)_{j+\frac12,k}^{n+1}
\right\rangle_h
=0\ ,\label{eq:modified-main-scheme1}
\\[1mm]
&
\delta_t^+ g_{j,k}^n
+
\frac1\varepsilon
(I-\Pi_h)\mathcal A_{j,k}^n
=
\frac1{\varepsilon^2}
\left(
K_{j,k}^{n,n+1}
-
g_{j,k}^{n+1}
\right)
-
\frac1\varepsilon
(I-\Pi_h)(\psi_{j,k}^n g_{j,k}^n)\ ,\label{eq:modified-main-scheme2}
\\[1mm]
&
\delta_t^+ S_{j+1/2}^n
=
 D_S(\delta_x^2S^{n+1})_{j+1/2}
-\alpha S_{j+1/2}^{n+1}
+\beta\rho_{j+1/2}^{n+1}\ ,\label{eq:modified-main-scheme3}
\\[1mm]
&
\delta_t^+N_{j+\frac12}^{n}
=
 D_N(\delta_x^2N^{n+1})_{j+1/2}
-\gamma\rho_{j+1/2}^{n+1}N_{j+1/2}^{n+1}\ .\label{eq:modified-main-scheme4}
\end{align}
\end{subequations}
Here the microscopic transport is treated explicitly by the upwind
approximation
\[
\mathcal A_{j,k}^n
=
v_k^+(\delta_x g)_{j-\frac12,k}^{n}
-
v_k^-(\delta_x g)_{j+\frac12,k}^{n},
\]
where
$
v_k^+=\max\{v_k,0\},
\
v_k^-=\max\{-v_k,0\}.
$
The local equilibrium part is defined by
\begin{equation}\label{eq:modified-K-two-upwind}
K_{j,k}^{n,n+1}
=
-\frac{v_k}{|V|_h}(\delta_x\rho)_j^{n+1}
-
\frac{\rho_j^{[S],n+1}}{|V|_h}\psi_{S,j,k}^n
-
\frac{\rho_j^{[N],n+1}}{|V|_h}\psi_{N,j,k}^n .
\end{equation}
The two upwind values $\rho_j^{[S],n+1}$ and $\rho_j^{[N],n+1}$ are chosen according to the effective velocities
\(u_{S,j}^n\) and \(u_{N,j}^n\),
\begin{equation}\label{eq: upwind densities}
\rho_j^{[S],n+1}
=
\begin{cases}
\rho_{j-\frac12}^{n+1}, & u_{S,j}^n\ge 0,\\
\rho_{j+\frac12}^{n+1}, & u_{S,j}^n<0,
\end{cases}
\quad 
\rho_j^{[N],n+1}
=
\begin{cases}
\rho_{j-\frac12}^{n+1}, & u_{N,j}^n\ge 0,\\
\rho_{j+\frac12}^{n+1}, & u_{N,j}^n<0.
\end{cases}
\end{equation}
This choice is consistent with the drift part of the limiting flux since
\[
-\frac{B_{S,j}^n}{|V|_h}\rho_j^{[S],n+1}
-
\frac{B_{N,j}^n}{|V|_h}\rho_j^{[N],n+1}
=
u_{S,j}^n\rho_j^{[S],n+1}
+
u_{N,j}^n\rho_j^{[N],n+1}.
\]
Taking the
discrete velocity average of the second equation in
\eqref{eq:modified-main-scheme} and using
$
\langle K_{j,k}^{n,n+1}\rangle_h=0
$
gives
\[
\langle g_{j,k}^{n+1}\rangle_h
=
\frac{\varepsilon^2}{\varepsilon^2+\Delta t}
\langle g_{j,k}^{n}\rangle_h .
\]
Thus \(\langle g_{j,k}^{n}\rangle_h=0\) is propagated in time if it holds
initially.

\paragraph{Decoupled solution strategy}
Following the micro-macro strategy in \cite{lemou2008new}, the scheme can
be solved in a decoupled way. We first introduce an auxiliary unknown
\(\widetilde g_{j,k}^{n+1}\) by
\begin{equation}\label{eq:modified-gtilda}
\frac{\widetilde g_{j,k}^{n+1}-g_{j,k}^n}{\Delta t}
+
\frac1\varepsilon
(I-\Pi_h)\mathcal A_{j,k}^n
=
\frac1{\varepsilon^2}
\left(
\widetilde K_{j,k}^{n}
-
\widetilde g_{j,k}^{n+1}
\right)
-
\frac1\varepsilon
(I-\Pi_h)(\psi_{j,k}^n g_{j,k}^n).
\end{equation}
where
\[
\widetilde K_{j,k}^{n}
=
-\frac{v_k}{|V|_h}(\delta_x\rho)_j^{n}
-
\frac{\rho_j^{[S],n}}{|V|_h}\psi_{S,j,k}^n
-
\frac{\rho_j^{[N],n}}{|V|_h}\psi_{N,j,k}^n .
\]
The upwind values \(\rho_j^{[S],n}\) and \(\rho_j^{[N],n}\) are defined in
the same way as in \eqref{eq: upwind densities}, with \(\rho^{n+1}\) replaced by \(\rho^n\).

All terms in \eqref{eq:modified-gtilda} are known except
\(\widetilde g_{j,k}^{n+1}\). Hence
\[
\widetilde g_{j,k}^{n+1}
=
\frac{\varepsilon^2}{\varepsilon^2+\Delta t}g_{j,k}^n
-
\frac{\varepsilon\Delta t}{\varepsilon^2+\Delta t}
(I-\Pi_h)\mathcal A_{j,k}^n
+
\frac{\Delta t}{\varepsilon^2+\Delta t}
\widetilde K_{j,k}^{n}
-
\frac{\varepsilon\Delta t}{\varepsilon^2+\Delta t}
(I-\Pi_h)(\psi_{j,k}^n g_{j,k}^n).
\]
Comparing \eqref{eq:modified-gtilda} with 
\eqref{eq:modified-main-scheme2}, we obtain
\begin{equation}\label{eq:modified-g-reconstruction}
g_{j,k}^{n+1}
=
\widetilde g_{j,k}^{n+1}
+
\lambda_\varepsilon
\left(
K_{j,k}^{n,n+1}
-
\widetilde K_{j,k}^{n}
\right),
\end{equation}
where
\begin{equation}\label{eq: lambda_eps}
\lambda_\varepsilon
=
\frac{\Delta t}{\varepsilon^2+\Delta t}.
\end{equation}

Substituting \eqref{eq:modified-g-reconstruction} into  \eqref{eq:modified-main-scheme1}, we obtain a closed equation for
\(\rho^{n+1}\). Define
\[
D_{\rho,h}
=
\frac{\langle v_k^2\rangle_h}{|V|_h},
\qquad
b_{S,j}^n
=
\frac{B_{S,j}^n}{|V|_h},
\qquad
b_{N,j}^n
=
\frac{B_{N,j}^n}{|V|_h},
\]
and let
$
a^n=\lambda_\varepsilon D_{\rho,h},
\
c_{S,j}^n=\lambda_\varepsilon b_{S,j}^n,
\
c_{N,j}^n=\lambda_\varepsilon b_{N,j}^n.
$
Then the density equation can be written as
\begin{equation}\label{eq:modified-rho-equation-two-upwind}
\frac{\rho_{j+\frac12}^{n+1}}{\Delta t}
-
\delta_x\left(
a^n\delta_x\rho^{n+1}
\right)_{j+\frac12}
-
\delta_x\left(
c_S^n\rho^{[S],n+1}
\right)_{j+\frac12}
-
\delta_x\left(
c_N^n\rho^{[N],n+1}
\right)_{j+\frac12}
=
r_{j+\frac12}^n,
\end{equation}
where the known right-hand side is
\[
\begin{aligned}
r_{j+\frac12}^n
=&
\frac{\rho_{j+\frac12}^n}{\Delta t}
-
\left\langle
v_k(\delta_x\widetilde g)_{j+\frac12,k}^{n+1}
\right\rangle_h
-
\delta_x\left(
a^n\delta_x\rho^n
\right)_{j+\frac12}
-
\delta_x\left(
c_S^n\rho^{[S],n}
\right)_{j+\frac12}
-
\delta_x\left(
c_N^n\rho^{[N],n}
\right)_{j+\frac12}.
\end{aligned}
\]
Thus only a scalar linear system for \(\rho^{n+1}\) needs to be solved at
each time step. With the upwind choice for the drift terms \eqref{eq: upwind densities}, the
resulting matrix is tridiagonal. Once \(\rho^{n+1}\) is obtained,
\(g^{n+1}\) is recovered explicitly from
\eqref{eq:modified-g-reconstruction}. 
%Finally, \(S^{n+1}\) and \(N^{n+1}\) are updated by  \eqref{eq:modified-main-scheme3} and \eqref{eq:modified-main-scheme4}.
We impose zero numerical density flux at both boundaries.
This condition is preserved by the positivity correction
described in Appendix~\ref{sec:limiter}.
The corrected density is then used to update $S^{n+1}$ and $N^{n+1}$
independently by \eqref{eq:modified-main-scheme3} and
\eqref{eq:modified-main-scheme4}, respectively.

\subsection{Asymptotic preserving property}

We now verify that the scheme is formally asymptotic preserving. 
Let \(\Delta t\) and \(\Delta x\) be fixed. Notice that $\lambda_\varepsilon$ in \eqref{eq: lambda_eps}
 converges to \(1\) as \(\varepsilon\) goes to zero.
Moreover, by the definition of \(\widetilde g^{n+1}\),
\(\widetilde g_{j,k}^{n+1}\) converges to \(\widetilde K_{j,k}^{n}\).
Using the expression of \(\widetilde K^n\), we obtain \(\left\langle
v_k(\delta_x\widetilde g)_{j+\frac12,k}^{n+1}
\right\rangle_h\) converges to
\[
-\delta_x\left(D_{\rho,h}\delta_x\rho^n\right)_{j+\frac12}
-
\delta_x\left(b_S^n\rho^{[S],n}\right)_{j+\frac12}
-
\delta_x\left(b_N^n\rho^{[N],n}\right)_{j+\frac12}.
\]
It follows from the definition of \(r_{j+\frac12}^n\) in
\eqref{eq:modified-rho-equation-two-upwind} that
\(
r_{j+\frac12}^n
\quad\text{converges to}\quad
\frac{\rho_{j+\frac12}^n}{\Delta t}.
\)
Therefore, in the limit as \(\varepsilon\) goes to zero,
\eqref{eq:modified-rho-equation-two-upwind} reduces to
\begin{equation}\label{eq:limiting-density-scheme}
\frac{\rho_{j+\frac12}^{n+1}-\rho_{j+\frac12}^{n}}{\Delta t}
-
\delta_x\left(
D_{\rho,h}\delta_x\rho^{n+1}
\right)_{j+\frac12}
-
\delta_x\left(
b_S^n\rho^{[S],n+1}
\right)_{j+\frac12}
-
\delta_x\left(
b_N^n\rho^{[N],n+1}
\right)_{j+\frac12}
=0 .
\end{equation}
This is a consistent implicit finite difference discretization of the
limiting drift-diffusion equation
\[
\partial_t\rho
-
\partial_x\left(D_{\rho,h}\partial_x\rho\right)
-
\partial_x\left(b_S\rho\right)
-
\partial_x\left(b_N\rho\right)
=0\ ,
\]
where \(u_S=-b_S\) and \(u_N=-b_N\).
Hence, for fixed $\Delta t$ and $\Delta x$, the fully discrete scheme reduces to a consistent discretization of the macroscopic drift–diffusion limit as $\varepsilon\to0$, which is the formal AP property in the diffusive regime.

\subsection{Bound-preserving property}

We now discuss the bound-preserving property of the limiting scheme.
In general, the full kinetic scheme is not automatically positivity
preserving. In the computation, this issue is handled by a post-processing
limiter as explained in Appendix \ref{sec:limiter}. For the limiting macroscopic scheme, however, a positivity and
boundedness result can be established.

In this subsection, we consider the limiting macroscopic scheme \eqref{eq:limiting-density-scheme} with
homogeneous no-flux boundary conditions, so that the discrete boundary
fluxes vanish. 

% Recall that the limiting scheme reads
% \[
% \delta_t^+\rho_{j+\frac12}^n
% =
% D_{\rho,h}\delta_x^2\rho_{j+\frac12}^{n+1}
% -
% \delta_x\left(u_S^n\rho^{[S],n+1}\right)_{j+\frac12}
% -
% \delta_x\left(u_N^n\rho^{[N],n+1}\right)_{j+\frac12},
% \]
% where the upwind values are chosen according to $u_{S,j}^n$ and $u_{N,j}^n$ as in \eqref{eq: upwind densities}

The limiting scheme can be written in matrix form as
\[
M^n\brho^{n+1}=\brho^n \ ,
\]
where
 $\boldsymbol{\rho}^{\,n}
=
(\rho_{j+\frac12}^n)_j$
denotes the vector of half-grid density values at time level $n$.
The tridiagonal matrix $M^n$ has entries
\begin{align*}
&m_{j,j-1}^n
=
-D_{\rho,h}\frac{\Delta t}{(\Delta x)^2}
-\frac{\Delta t}{\Delta x}(u_{S,j}^n)^+
-\frac{\Delta t}{\Delta x}(u_{N,j}^n)^+,\\
&m_{j,j}^n
=
1+2D_{\rho,h}\frac{\Delta t}{(\Delta x)^2}
+\frac{\Delta t}{\Delta x}
\left[(u_{S,j+1}^n)^+ +(u_{S,j}^n)^-\right]
+\frac{\Delta t}{\Delta x}
\left[(u_{N,j+1}^n)^+ +(u_{N,j}^n)^-\right],\\
&m_{j,j+1}^n
=
-D_{\rho,h}\frac{\Delta t}{(\Delta x)^2}
-\frac{\Delta t}{\Delta x}(u_{S,j+1}^n)^-
-\frac{\Delta t}{\Delta x}(u_{N,j+1}^n)^- .
\end{align*}
Here $(a)^+=\max\{a,0\}$ and $(a)^-=\max\{-a,0\}$.
The boundary rows of \(M^n\) are modified according to the homogeneous
no-flux boundary conditions.

\begin{proposition}[Mass conservation, positivity, and upper bound]
Assume that $D_{\rho,h}>0$ and that the effective velocities satisfy
\[
|u_{S,j}^n|\le U_S,\qquad |u_{N,j}^n|\le U_N
\]
for all $j$ and $n$. Suppose further that
\begin{equation}\label{eq:CFL-bound-new}
\frac{\Delta t}{\Delta x}
\le
\frac{1}{4(U_S+U_N)} .
\end{equation}
If the initial data satisfy $\rho_{j+\frac12}^0\ge0$ for all $j$, then the
limiting scheme satisfies
\[
\rho_{j+\frac12}^n\ge0,\qquad \forall j,\ n\ge0.
\]
Moreover, under the homogeneous no-flux boundary conditions, the discrete
mass is conserved,
\[
\|\brho^{n+1}\|_1=\|\brho^n\|_1 := \sum_j|\rho^n_{j+\frac12}|,
\]
and the solution satisfies the uniform bound in time on a fixed spatial grid
\[
0\le \rho_{j+\frac12}^n\le \bar\rho,\qquad
\bar\rho=2\|\brho^0\|_1 .
\]
\end{proposition}

\begin{proof}{\rm
We first observe that
\[
m_{j,j-1}^n+m_{j,j}^n+m_{j,j+1}^n
=
1+\frac{\Delta t}{\Delta x}
\left(u_{S,j+1}^n-u_{S,j}^n\right)
+\frac{\Delta t}{\Delta x}
\left(u_{N,j+1}^n-u_{N,j}^n\right)
=:F_j^n .
\]
Using the bounds on $u_S^n$ and $u_N^n$, we have
\[
F_j^n
\ge
1-2(U_S+U_N)\frac{\Delta t}{\Delta x}.
\]
Under \eqref{eq:CFL-bound-new}, this gives
$
F_j^n\ge \frac12 .
$

The off-diagonal entries of $M^n$ are nonpositive. Moreover, the above
lower bound implies the required diagonal dominance under the CFL
condition. Hence $M^n$ is a nonsingular M-matrix and $(M^n)^{-1}$ is
nonnegative. Therefore, if $\brho^n\ge0$, then
\[
\brho^{n+1}=(M^n)^{-1}\brho^n\ge0.
\]
Mass conservation follows from the conservative form of the scheme and the
homogeneous no-flux boundary conditions.  More precisely, 
\[
\|\brho^{n+1}\|_1
=
\mathbf 1^T\brho^{n+1}
=
\mathbf 1^T M^n\brho^{n+1}
=
\mathbf 1^T\brho^n
=
\|\brho^n\|_1 ,
\]
where we used the relation
\(
\mathbf 1^T M^n=\mathbf 1^T .
\)

Finally, let $\bar\rho=2\|\brho^0\|_1$. By mass conservation and
positivity,
\[
\rho_{j+\frac12}^n\le \|\brho^n\|_1=\|\brho^0\|_1=\frac12\bar\rho .
\]
Since $F_j^n\ge1/2$, we obtain
\[
(M^n(\bar\rho\mathbf 1))_j
=
F_j^n\bar\rho
\ge
\frac12\bar\rho
\ge
\rho_{j+\frac12}^n .
\]
Using the nonnegativity of $(M^n)^{-1}$, we conclude that
\[
\brho^{n+1}
=
(M^n)^{-1}\brho^n
\le
(M^n)^{-1}M^n(\bar\rho\mathbf 1)
=
\bar\rho\mathbf 1 .
\]
This proves the desired upper bound. }
\end{proof}

\section{Numerical experiments}\label{sec: numerics}

In this section, we present numerical experiments to illustrate the
traveling pulse behavior of the model and to verify the asymptotic
preserving property of the proposed micro-macro scheme. We first validate
the macroscopic solver by comparing the numerical traveling pulse with
the analytical traveling profile. We then study the dependence of the
traveling speed on the chemotactic sensitivities. Finally, we compare
the kinetic and macroscopic solutions for different values of
\(\varepsilon\).

Unless otherwise stated, we consider the one-dimensional computational
domain \(\Omega=(0,150)\). For the profile and asymptotic preserving
tests, we use
\(
\Delta x=0.025
\), 
\(
\Delta t=10^{-3}.
\)
The velocity domain is \(V=[-1,1]\). It is discretized by a uniform
trapezoidal quadrature with \(\Delta v=0.02\). We denote the discrete
velocity nodes and weights by \(\{v_k,\omega_k\}\), and use
\[
\langle h\rangle_h=\sum_k\omega_k h(v_k),
\qquad
|V|_h=\sum_k\omega_k .
\]
The discrete diffusion coefficient in the limiting macroscopic equation
is
 $D_{\rho,h}
=
\frac{\langle v^2\rangle_h}{|V|_h}.$
For the default velocity grid with \(\Delta v=0.02\), this gives
\(D_{\rho,h}=0.3334\).

We use the reflecting boundary conditions described in
Appendix~\ref{sec:limiter} for the kinetic variable and
homogeneous Neumann conditions for the chemical variables.

The initial data are chosen so that cells are concentrated near the left
side of the domain, with no chemoattractant signal and a constant
nutrient level. More precisely, we take
\[
\rho(0,x)=\rho_0 e^{-\lambda x},
\qquad
S(0,x)=0,
\qquad
N(0,x)=N_0,
\]
where
$\rho_0=1,
\
\lambda=3,
\
N_0=10^3 .$
For the kinetic computations, the initial data are chosen consistently
with the micro-macro decomposition,
\[
f^\varepsilon(0,x,v)
=
\frac{\rho(0,x)}{|V|_h}
+
\varepsilon g(0,x,v),
\qquad
g(0,x,v)=0 .
\]
The chemical parameters are
$
D_S=2,
\
D_N=0,
\
\alpha=0.05,
\
\beta=\gamma=1 .
$
Unless otherwise stated, the chemotactic sensitivities in the profile
and AP tests are
$\chi_S=0.5,
\
\chi_N=1.1 .$
For the experiments in
Sections~\ref{sec:macro-validation}--\ref{sec:sign-ap-experiments},
we use
\[
\phi_h(y)=-s_h\operatorname{sign}(y),
\qquad
s_h=\frac{|V|_h}{\langle|v|\rangle_h}.
\]
With this normalization, the macroscopic drift velocities are
$u_S=\chi_S\operatorname{sign}(S_x)$ and
$u_N=\chi_N\operatorname{sign}(N_x)$.
The response arguments are evaluated with $\varepsilon=0$
for the macroscopic model and with $\varepsilon>0$ for the kinetic model.
The kinetic computations use the correction described in Appendix~\ref{sec:limiter}.

\subsection{Macroscopic traveling pulse and analytical profile}
\label{sec:macro-validation}

We first validate the macroscopic solver by comparing the numerical
traveling pulse with the analytical traveling profile. Let
\(x_p\) be the peak position of the numerical density and set
\[
z=x-x_p,
\qquad
z_-=\min\{z,0\},
\qquad
z_+=\max\{z,0\}.
\]
The analytical traveling profile with the response function $\phi(y)=-\operatorname{sign}(y)$ found in \cite{calvez2019chemotactic} is written as 
\begin{equation}
\rho_{\rm tw}(x)
=
\rho_p
\exp\left(
\lambda_L z_-
+
\lambda_R z_+
\right),
\label{eq:analytic-sol}
\end{equation}
where
\[
\lambda_L
=
\frac{\chi_S+\chi_N-\sigma}{D_{\rho,h}},
\qquad
\lambda_R
=
\frac{\chi_N-\chi_S-\sigma}{D_{\rho,h}} .
\]
Here \(\lambda_L>0\) and \(\lambda_R<0\), so that the profile decays on
both sides of the peak. The amplitude \(\rho_p\) is chosen by discrete
mass matching. More precisely, if \(z_j=x_{j+\frac12}-x_p\), then
\[
\rho_p
=
\frac{
\Delta x\sum_j \rho_j
}{
\Delta x\sum_j
\exp\left(
\lambda_L (z_j)_-
+
\lambda_R (z_j)_+
\right)
}.
\]

Figure~\ref{fig:macro_analytic_cmp} compares the macroscopic numerical
profiles with the analytical traveling profiles at two representative
times. The peak position and amplitude of the analytical profile are
chosen from the corresponding numerical solution. The close agreement
confirms that the macroscopic solver captures the expected traveling
wave structure. This macroscopic solution will be used as the reference
solution in the kinetic-to-macroscopic comparison below.

\begin{figure}[htbp]
    \centering
    \begin{subfigure}{0.48\linewidth}
        \centering
        \includegraphics[width=\linewidth]{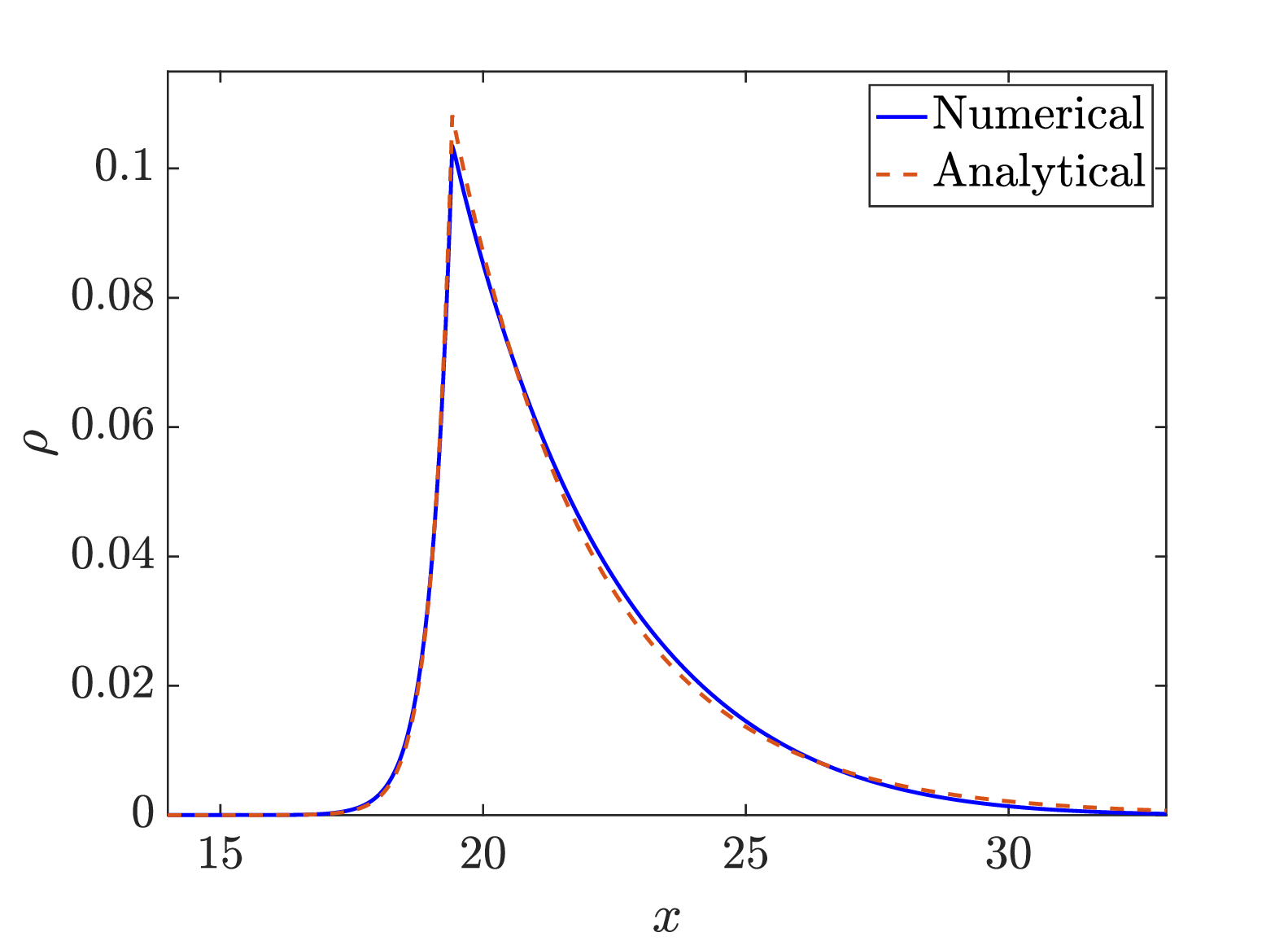}
        \caption{\(t=30\)}
        \label{fig:macro_analytic_t30}
    \end{subfigure}
    \begin{subfigure}{0.48\linewidth}
        \centering
        \includegraphics[width=\linewidth]{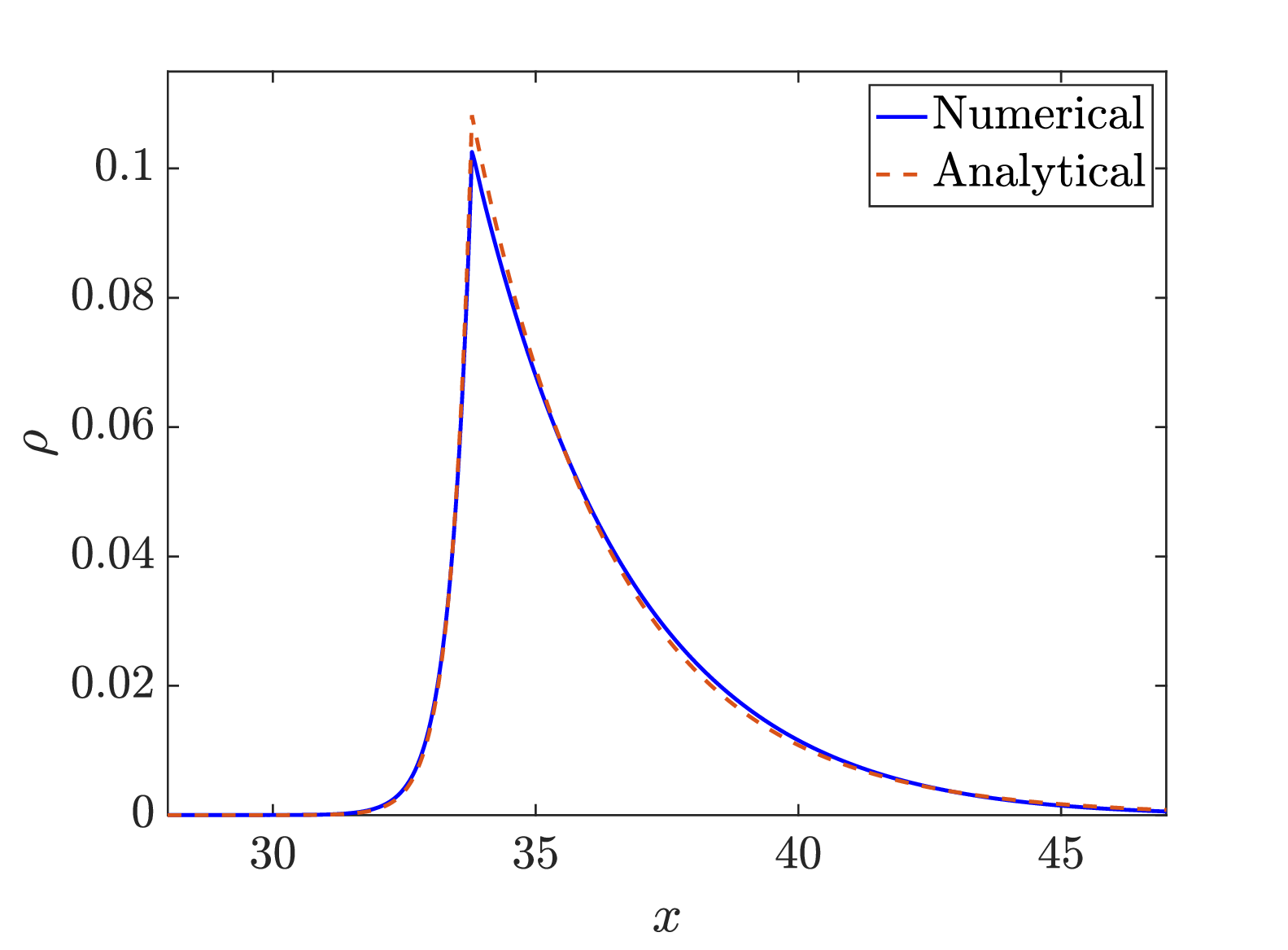}
        \caption{\(t=50\)}
        \label{fig:macro_analytic_t50}
    \end{subfigure}
    \caption{Comparison between the macroscopic numerical profiles and
    the analytical traveling profiles at \(t=30\) and \(t=50\). }
    \label{fig:macro_analytic_cmp}
\end{figure}

\subsection{Dependence of the traveling speed on chemotactic sensitivities}

We next investigate how the traveling speed depends on \(\chi_S\) and
\(\chi_N\). For each pair \((\chi_S,\chi_N)\), the predicted traveling
speed \(\sigma\) is determined by
\begin{equation}
\chi_N-\sigma
=
\chi_S
\frac{\sigma}{\sqrt{4D_S\alpha+\sigma^2}} .
\label{eq:sigma-macro}
\end{equation}
To estimate the traveling speed numerically, we track the center of
mass of the density,
\begin{equation}\label{eq: center of mass}
X_\rho(t_n)
=
\frac{
\Delta x\sum_j x_{j+\frac12}\rho_{j+\frac12}^n
}{
\Delta x\sum_j \rho_{j+\frac12}^n
}.
\end{equation}

For the parameter study, we use a slightly coarser
mesh to reduce the computational cost, i.e. 
\(
\Delta x=0.1,
\Delta v=0.02,
\Delta t=10^{-2}.
\)
The sensitivities \(\chi_S,\chi_N\) are uniformly sampled on \([0,2] \times [0,2]\). 
The computational domain for this parameter study is $\Omega=(0,300)$. 
The wave speeds are evaluated at $t=100$.  

Figure~\ref{fig:traveling_speed} shows the dependence of the traveling
speed on \(\chi_S\) and \(\chi_N\). The subfigure (a) gives the speed
predicted by \eqref{eq:sigma-macro}, while the subfigures (b) and (c) show the
numerical speed obtained from the macroscopic and the kinetic model with $\varepsilon=0.01$, respectively, by tracking the center of mass. 
The three panels agree well in the parameter regime where
a clear traveling pulse is formed. Some discrepancies are observed near
the transition region where the nutrient bias is not sufficiently strong
compared with the self-attractant bias. In this region, the pulse
formation may require a longer transient time, and the finite-time
center-of-mass estimate may underestimate the asymptotic traveling
speed. 
%(XR: Yes, the plot for the kinetic case is almost the same. I removed it because the purpose of this subsection is mainly to validate the macroscopic speed formula \eqref{eq:sigma-macro}.The kinetic speed comparison is now included in Figure~\ref{fig:conv_speed}, where we show its convergence to the macroscopic speed as \(\varepsilon\) decreases.)

\begin{figure}[htbp]
    \centering
    \begin{subfigure}{0.3\linewidth}
        \centering
        \includegraphics[width=\linewidth]{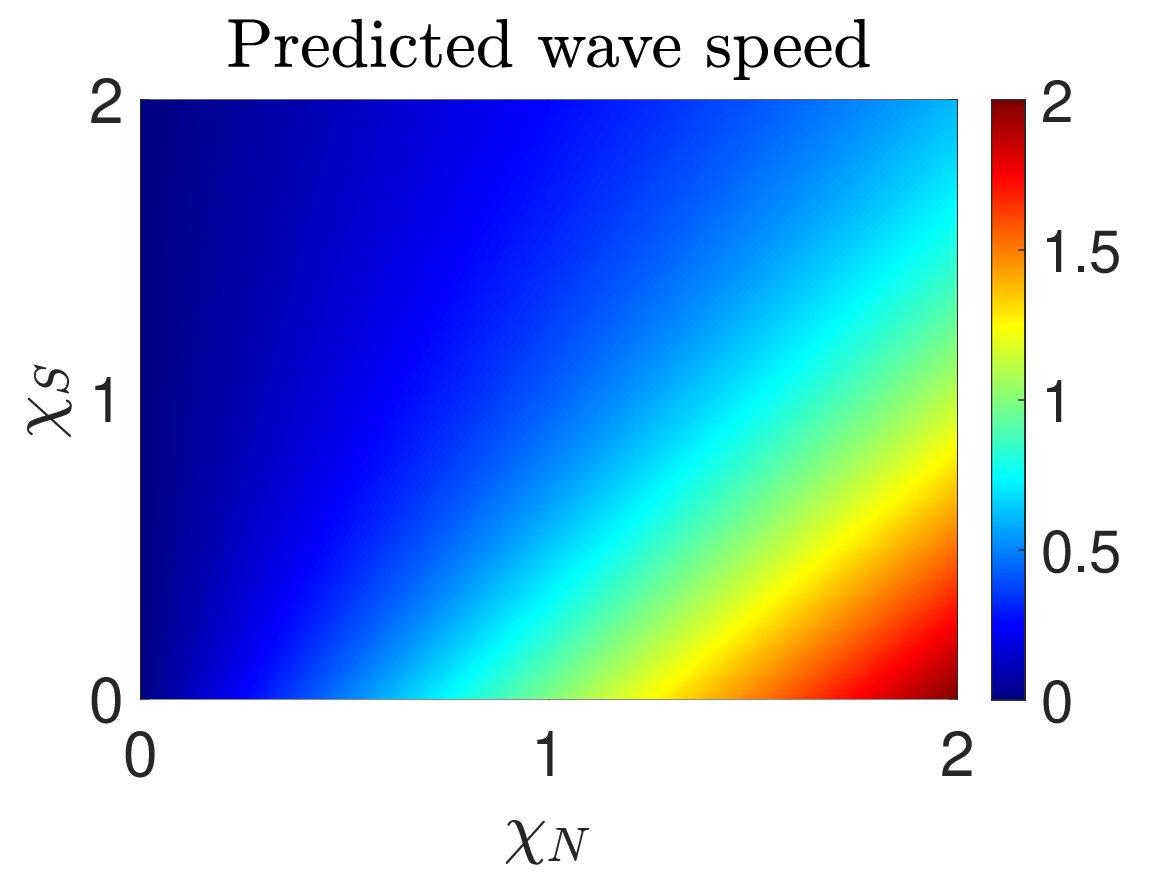}
        \caption{}
        \label{fig:analytic_sigma}
    \end{subfigure}
    \begin{subfigure}{0.3\linewidth}
        \centering
        \includegraphics[width=\linewidth]{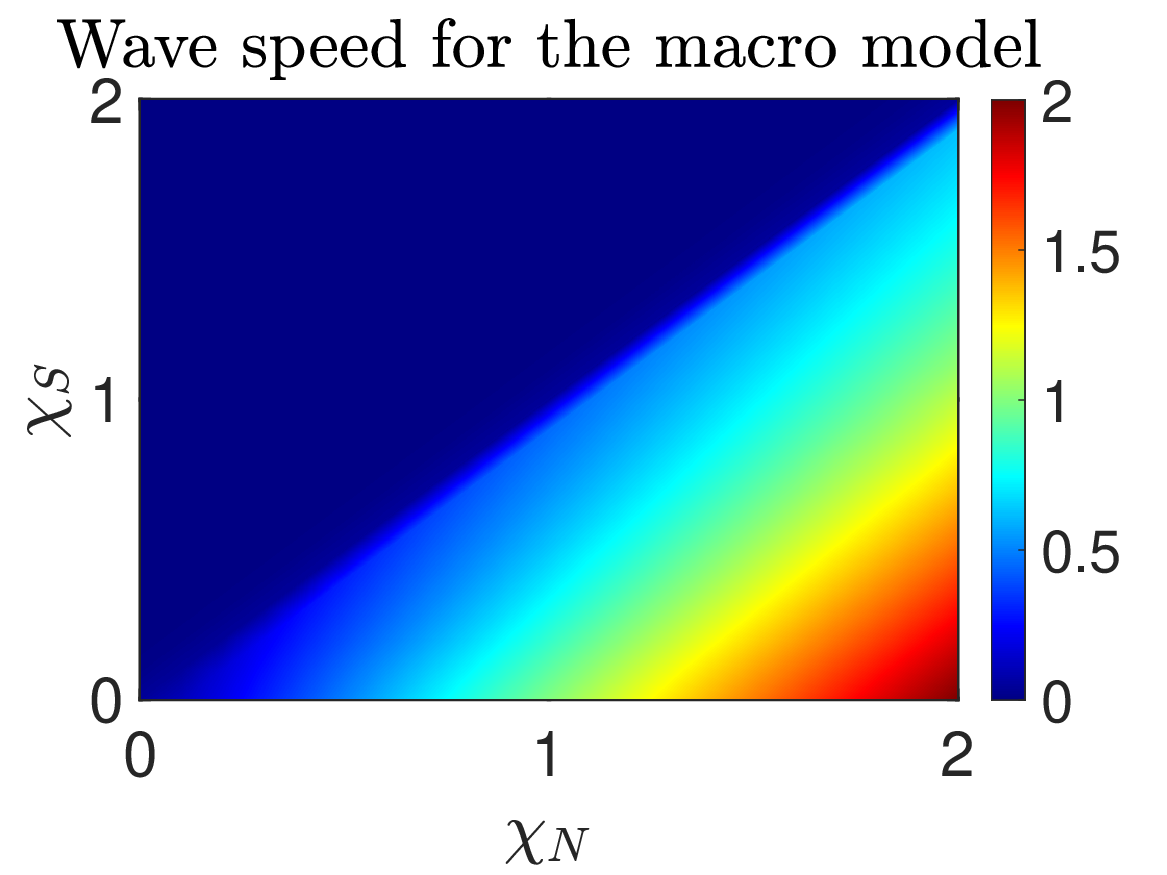}
        \caption{}
        \label{fig:numerical_macro_sigma}
    \end{subfigure}
    \begin{subfigure}{0.3\linewidth}
        \centering
        \includegraphics[width=\linewidth]{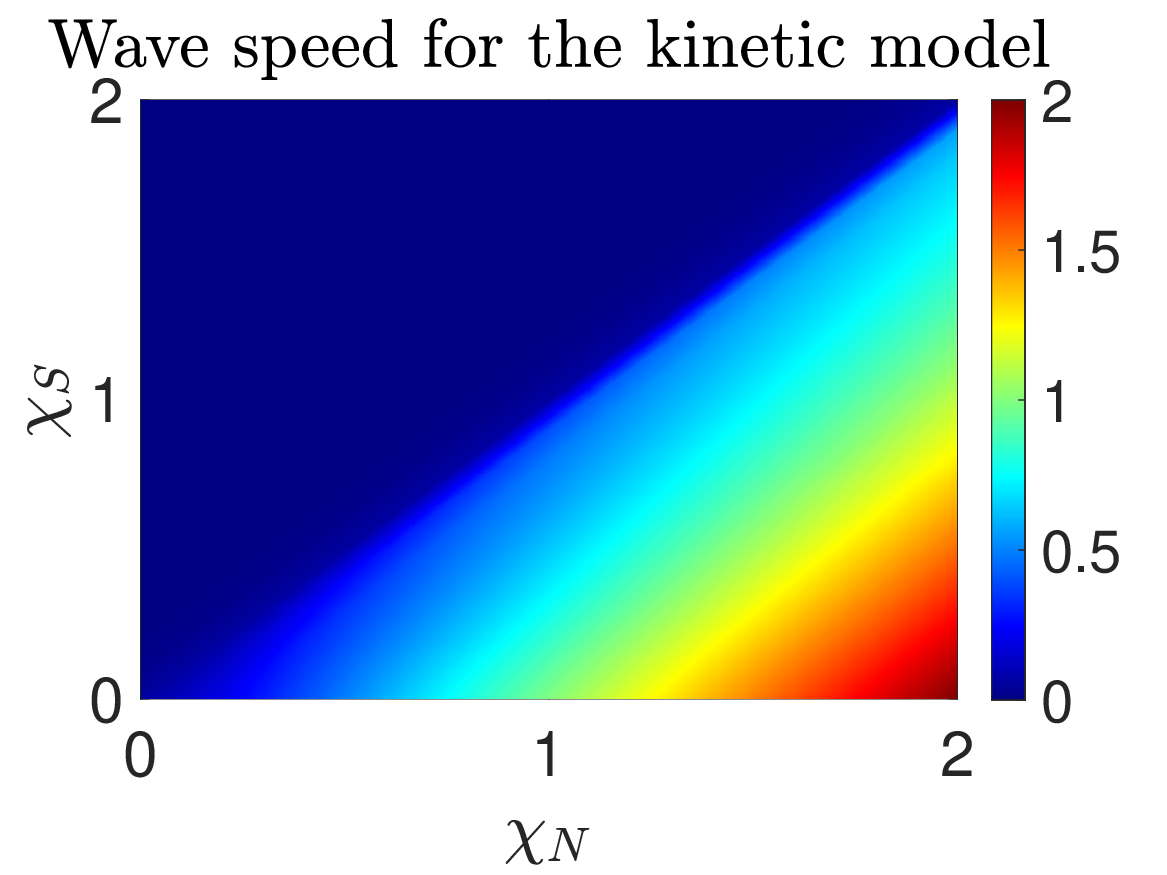}
        \caption{}
        \label{fig:numerical_kinetic_sigma}
    \end{subfigure}
    \caption{Dependence of the traveling speed on \(\chi_S\) and
    \(\chi_N\). Figure (a) shows the speed predicted by
    \eqref{eq:sigma-macro}, while Figure (b) and Figure (c) show the numerical
    speed obtained from the macroscopic and kinetic model with $\varepsilon=0.01$, respectively.}
    \label{fig:traveling_speed}
\end{figure}

\subsection{Asymptotic preserving convergence} \label{sec:sign-ap-experiments}

We now compare the kinetic solution with the limiting macroscopic
solution for different values of \(\varepsilon\). The purpose of this
test is to verify that, for fixed \(\Delta x\), \(\Delta v\), and
\(\Delta t\), the kinetic AP solution approaches the limiting
macroscopic solution as \(\varepsilon\) decreases. We use the baseline
parameter setting described at the beginning of this section and consider
$\varepsilon=0.1,
\
\varepsilon=0.05,
\
\varepsilon=0.01.$

Figure~\ref{fig:conv_T} compares the traveling density pulses obtained
from the kinetic model and the limiting macroscopic model at several
representative times. As \(\varepsilon\) decreases, the kinetic profiles
become progressively closer to the macroscopic profile. In particular,
for \(\varepsilon=0.01\), the kinetic and
macroscopic solutions are already nearly indistinguishable at the
plotted times. This provides direct evidence that the AP scheme captures
the correct diffusive limit.

\begin{figure}[htbp]
    \centering
    \includegraphics[width=\linewidth]{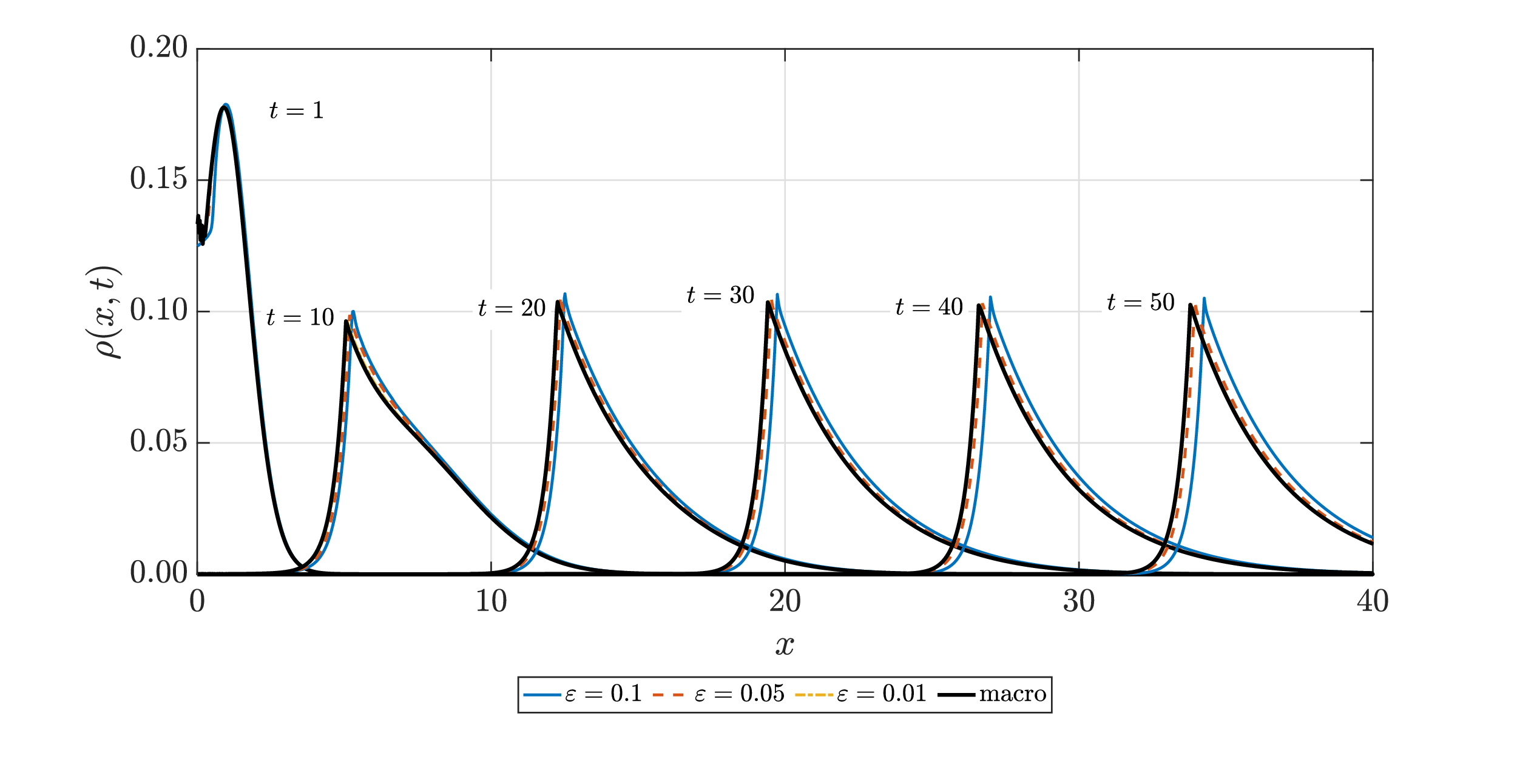}
    \caption{Comparison of traveling density pulses obtained from the
    kinetic model and the limiting macroscopic model at
    \(t=1,10,20,30,40,50\). The propagation of the pulse is captured at
    both levels, and the kinetic profiles approach the macroscopic
    profiles as \(\varepsilon\) decreases. }
    \label{fig:conv_T}
\end{figure}

To quantify the convergence, we use the discrete \(L^2\)-error
\[
\operatorname{err}_{\rho,\varepsilon}(t)
=
\left(
\Delta x\sum_j
\left|
\rho_{j+\frac12}^{\varepsilon}(t)
-
\rho_{j+\frac12}(t)
\right|^2
\right)^{1/2},
\]
where \(\rho^\varepsilon\) is reconstructed from the kinetic solution
and \(\rho\) denotes the solution of the limiting macroscopic model.
Figure~\ref{fig:conv_error} plots the final-time error
\(\operatorname{err}_{\rho,\varepsilon}(50)\) for different values of
\(\varepsilon\). The error decreases monotonically as \(\varepsilon\)
becomes smaller, confirming the convergence of the kinetic AP solution
toward the macroscopic one.

\begin{figure}[htbp]
    \centering
    \includegraphics[width=0.6\linewidth]{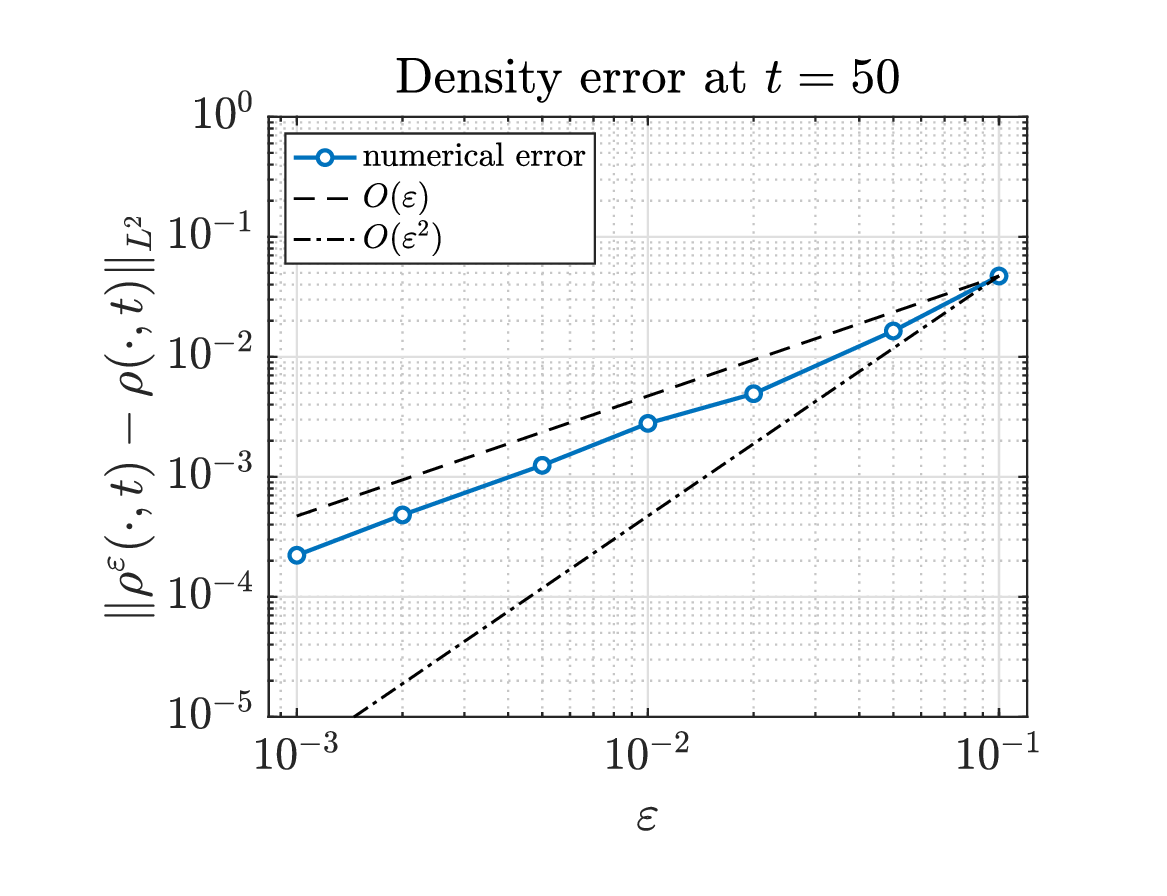}
    \caption{Final-time density error
    \(\operatorname{err}_{\rho,\varepsilon}(50)\) versus
    \(\varepsilon\). The reference lines \(O(\varepsilon)\) and
    \(O(\varepsilon^2)\) are included for comparison.}
    \label{fig:conv_error}
\end{figure}

Besides the density profile, it is also important to compare the
traveling speed. The instantaneous speed is computed from the
center-of-mass trajectory \(X_\rho(t)\) introduced in \eqref{eq: center of mass}. Figure~\ref{fig:conv_speed}
shows the speeds obtained from the macroscopic solution and from the
kinetic solutions for different values of \(\varepsilon\). The horizontal
dashed line indicates the predicted value
$
\sigma_{\rm pred}=0.724,$
obtained from \eqref{eq:sigma-macro}. 
After a short transient regime, the propagation speeds
become approximately constant and remain close to the
predicted value. 
Moreover, the kinetic wave speed
becomes closer to the macroscopic one as \(\varepsilon\) decreases.

\begin{figure}[htbp]
    \centering
    \includegraphics[width=0.82\linewidth]{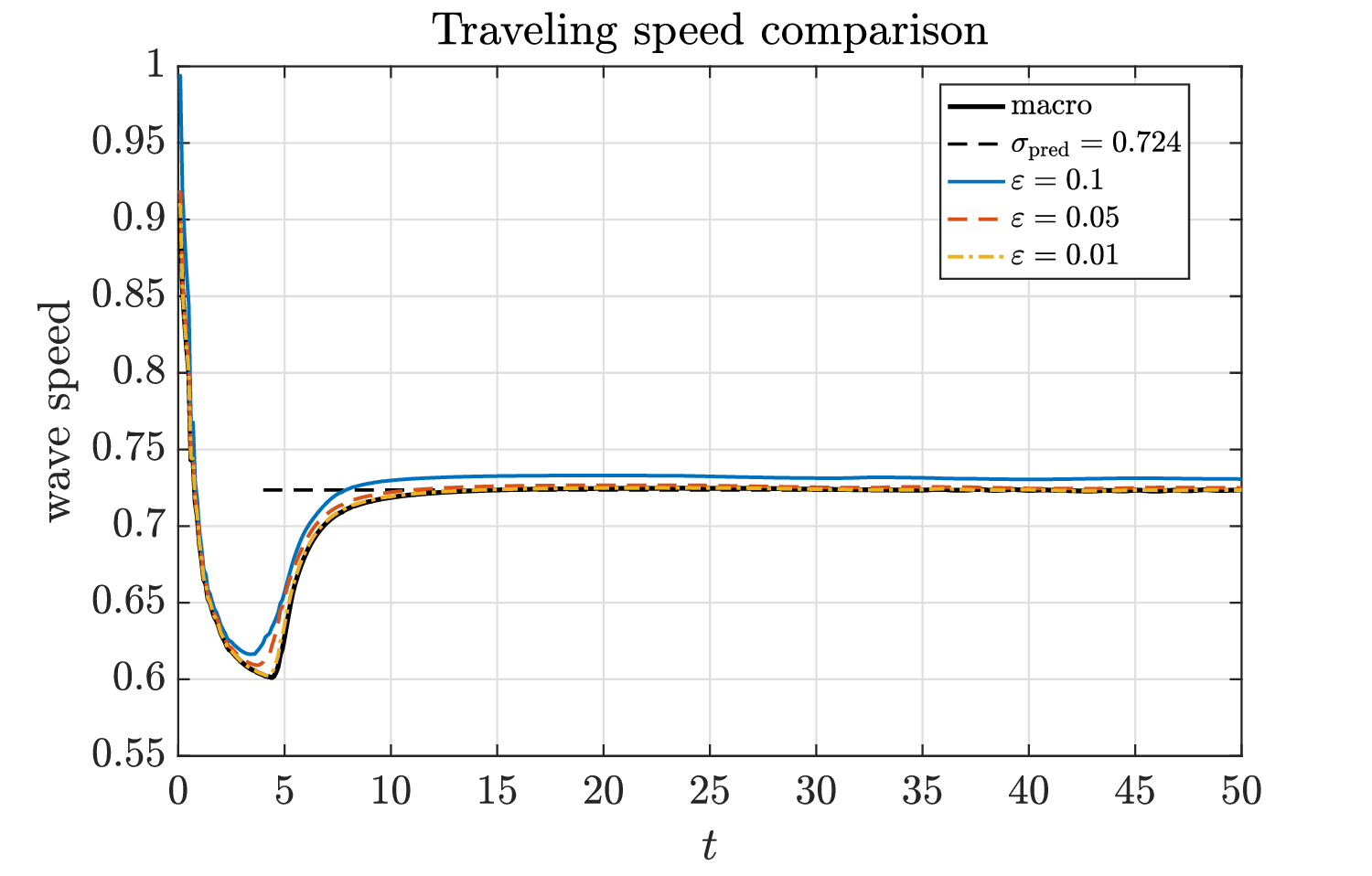}
    \caption{traveling speed comparison between the macroscopic model
    and the kinetic model for different values of \(\varepsilon\). 
    The horizontal dashed line indicates the macroscopic predicted speed \(\sigma_{\rm macro}\approx0.724\), obtained from \eqref{eq:sigma-macro}.}
    \label{fig:conv_speed}
\end{figure}

\subsection{Traveling pulses with an exponential response}
\label{sec:exponential-response}

We next consider the exponentially modulated response
\begin{equation}
\phi_{a,h}(y)
=
-s_h e^{-a|y|}\operatorname{sign}(y),
\qquad
s_h=\frac{|V|_h}{\langle|v|\rangle_h}.
\label{eq:exponential-response}
\end{equation}
Here $a\ge0$ controls the exponential modulation, and $a=0$
recovers the normalized sign response.
The initial data and all other parameters are the same as in
Section~\ref{sec:sign-ap-experiments}.
For each value of $a$, the macroscopic reference solution is
computed using the same response function with $\varepsilon=0$.

For $a=0.001$ and $\varepsilon=0.01$,
Figure~\ref{fig:exp-response-evolution} shows an asymmetric pulse
propagating to the right.
Its propagation speed becomes approximately constant at late times.
At fixed $a=0.001$, the kinetic density profiles approach the
corresponding numerical macroscopic solution as $\varepsilon$
decreases, as shown in Figure~\ref{fig:exp-response-limits}(a).
At fixed $\varepsilon=0.01$, the density profiles approach
the sign-response solution as $a$ decreases
in Figure~\ref{fig:exp-response-limits}(b).
The fitted propagation speeds also approach the numerical
sign-response value, as reported in
Table~\ref{tab:exp-response-speeds}.

\begin{figure}[htbp]
    \centering
    \color{red}
    % Replace with the actual EPS filename.
    \includegraphics[width=0.82\linewidth]{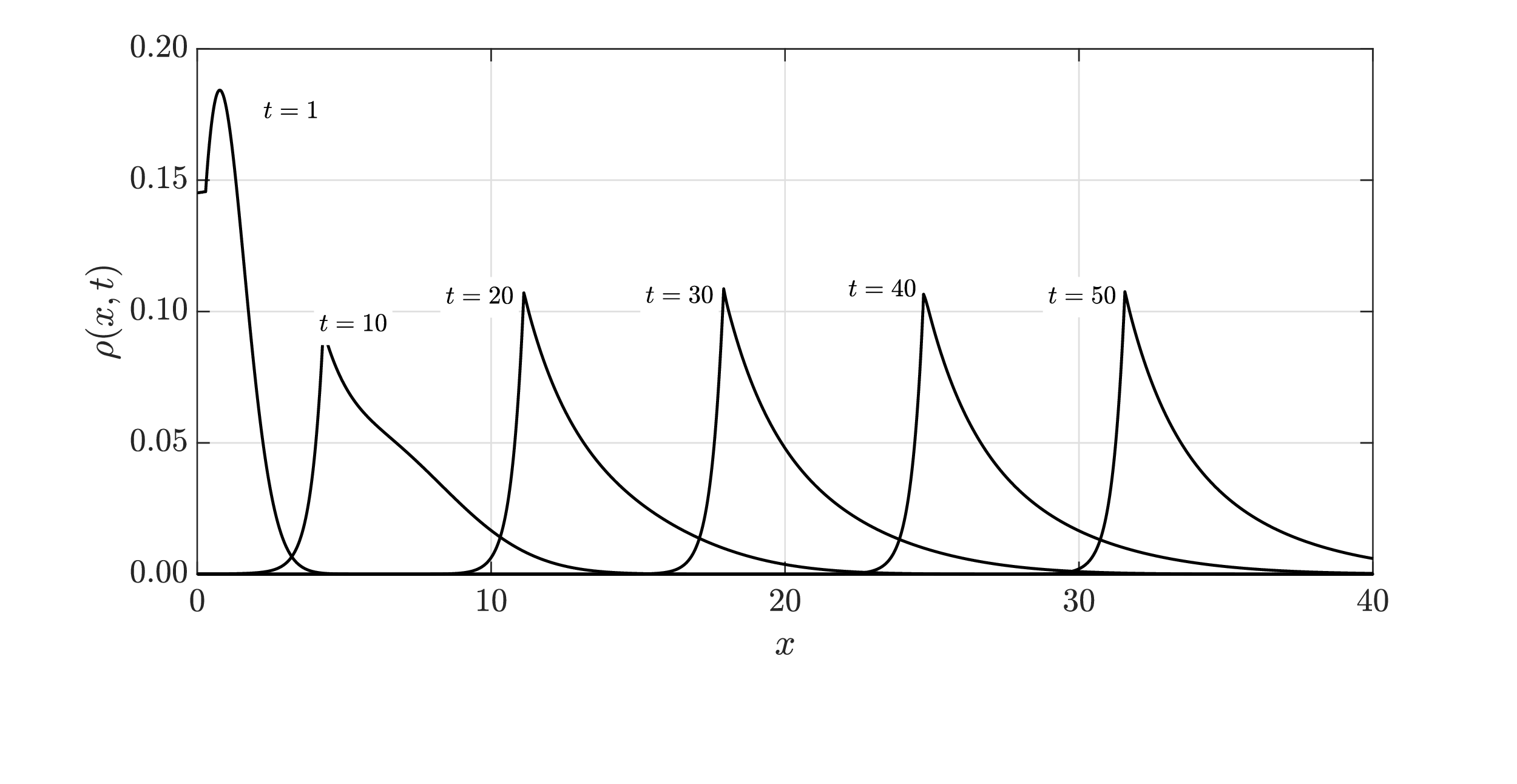}
    \caption{Density profiles for the exponential response
    with $a=0.001$ and $\varepsilon=0.01$
    at $t=1,10,20,30,40,50$.}
    \label{fig:exp-response-evolution}
\end{figure}

\begin{figure}[htbp]
    \centering
    \color{red}
    % Replace with the actual EPS filename.
    % This file contains both panels (a) and (b).
    \includegraphics[width=0.95\linewidth]{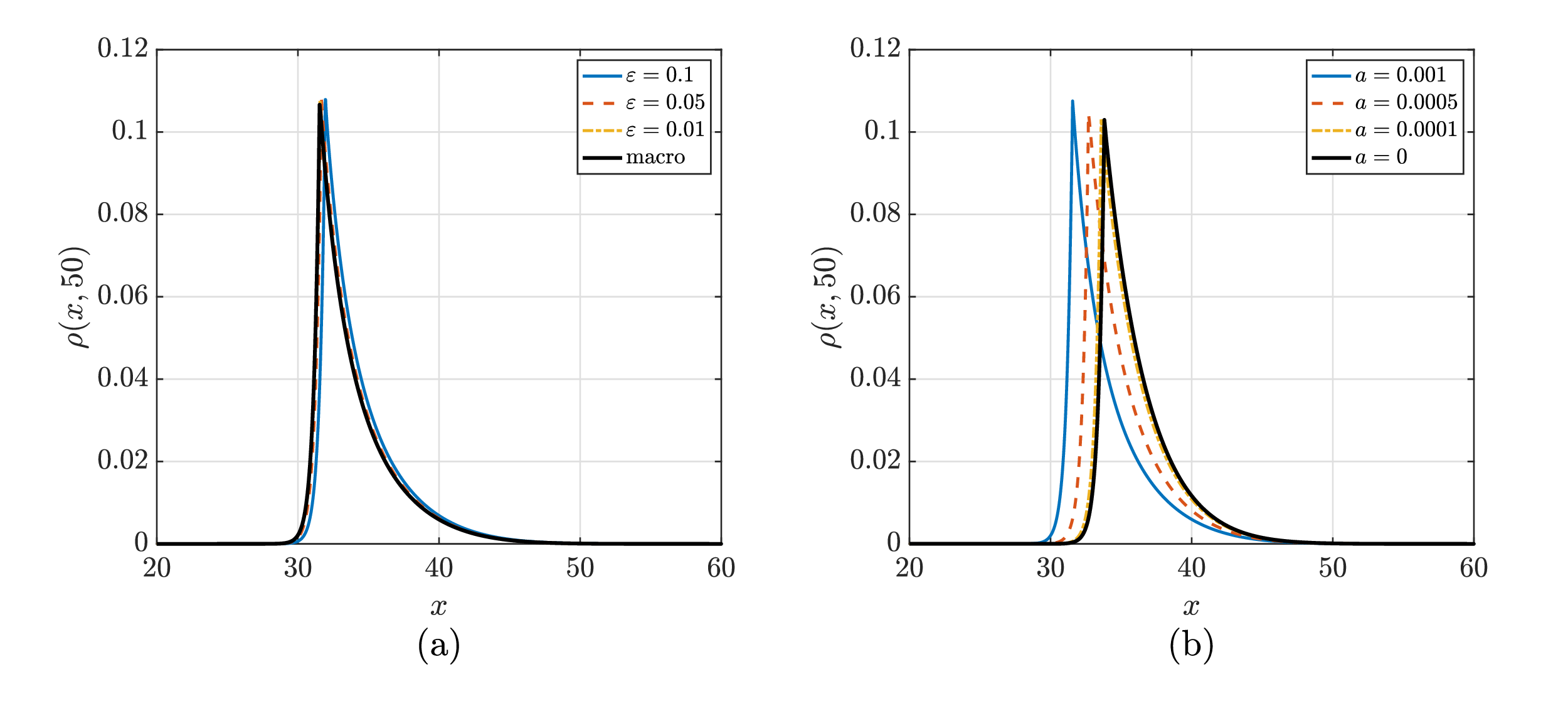}
    \caption{Density profiles at $T=50$.
    (a) Kinetic solutions with $\varepsilon=0.1,0.05,0.01$
    and the corresponding numerical macroscopic solution
    for $a=0.001$.
    (b) Solutions with $a=0.001,0.0005,0.0001,0$
    for $\varepsilon=0.01$.
    The case $a=0$ corresponds to the sign response.}
    \label{fig:exp-response-limits}
\end{figure}

\begin{table}[htbp]
    \centering
    \caption{Numerical propagation speeds for the exponential
    response with $\varepsilon=0.01$.
    The speeds are obtained by least-squares linear fits to the center of mass over $t\in[30,50]$.}
    \label{tab:exp-response-speeds}
    \begin{tabular}{lcccc}
        \hline
        $a$
        & $10^{-3}$ & $5\times10^{-4}$ & $10^{-4}$ & $0$ \\
        \hline
        Propagation speed
        & $0.690$ & $0.707$ & $0.720$ & $0.724$ \\
        \hline
    \end{tabular}
\end{table}

Overall, these results show that the proposed AP scheme captures the
traveling pulse profiles at the kinetic level and converges to the
limiting macroscopic dynamics, both in density profiles and in
propagation speed, as \(\varepsilon\) decreases.

\FloatBarrier

}

\section{Conclusions}\label{sec: conclusions}

In this work, we developed and analyzed an asymptotic preserving finite difference scheme for a kinetic chemotaxis model describing bacterial traveling pulses. The model couples a run-and-tumble kinetic equation with reaction--diffusion equations for a nutrient and a self-produced chemoattractant, motivated by experimental observations of \textit{E.~coli} pulse propagation in microchannels.

The proposed numerical method is based on a micro--macro decomposition of the kinetic equation and remains stable in both the kinetic and diffusive regimes. In the limit as the scaling parameter $\varepsilon \to 0$, the scheme naturally reduces to a consistent discretization of the limiting drift--diffusion system. We also established positivity, mass conservation, and uniform boundedness properties for the limiting macroscopic scheme, and introduced a positivity-preserving post-processing procedure for the kinetic solver.

The numerical experiments confirm that the scheme accurately captures the formation and propagation of traveling bacterial pulses for the sign response function. In particular, the simulations reproduce the characteristic asymmetric pulse profile and validate the theoretical prediction for the traveling wave speed found in \cite{calvez2019chemotactic}. Moreover, the kinetic solutions converge toward the macroscopic limit as $\varepsilon$ decreases, demonstrating the asymptotic preserving character of the method.

Several directions remain open for future work. A rigorous convergence analysis of the fully discrete AP scheme would be of significant interest. It would also be natural to extend the method to higher-dimensional settings and more general kinetic response functions. Another important direction concerns the development of fully positivity-preserving AP schemes without post-processing corrections. Finally, incorporating additional biological mechanisms, such as population growth, signal adaptation, or interactions between multiple bacterial species, could provide a more realistic description of collective bacterial migration phenomena.

\bigskip
\noindent{\large\bfseries Acknowledgements} 
CA was supported by TÜBİTAK (The Scientific and Technological Research Council of Türkiye) through the International Postdoctoral Research Fellowship Program under Grant No: 1059B191800629.
JAC was supported by the Advanced Grant Nonlocal-CPD (Nonlocal PDEs for Complex Particle Dynamics: Phase Transitions, Patterns and Synchronization) of the European Research Council Executive Agency (ERC) under the European Union Horizon 2020 research and innovation programme (grant agreement No. 883363) and  partially supported by the EPSRC EP/V051121/1. JAC was partially supported by the ``Maria de Maeztu'' Excellence Unit IMAG, reference CEX2020-001105-M, funded by MCIN/AEI/10.13039/501100011033/.
GER is supported by the Spanish grant PID2022-143012NA-100 and she is a member of the Catalan research group 2021-SGR-00087.
XR acknowledges support from the R\&D Program of Beijing Municipal Education Commission from China under Grant No.  KM202310028016, and the National Natural Science Foundation of China under Grant No. 12201436.

\appendix

\section{Boundary treatment and positivity preservation}
\label{sec:limiter}
We use the discrete velocity average
\[
\langle h\rangle_h=\sum_{k}\omega_k h(v_k),
\qquad \omega_k>0,\qquad
|V|_h=\sum_{k}\omega_k,\qquad
\psi_0(v_k)=|V|_h^{-1}.
\]
Let $J$ be the number of spatial cells, and assume
$v_{-k}=-v_k$ and $\omega_{-k}=\omega_k$.
We impose reflection at the endpoints and use reflected ghost values,
\[
g_{j,k}^n=g_{j,-k}^n\quad(j=0,J),\qquad
g_{-1,k}^n=g_{1,-k}^n,\qquad
g_{J+1,k}^n=g_{J-1,-k}^n.
\]
Endpoint symmetry is imposed before the kinetic update, on the
auxiliary variable $\widetilde g^{n+1}$, and on the provisional
microscopic variable before checking positivity.
The diffusion and drift boundary fluxes are set to zero.
Thus the complete numerical density flux vanishes at both endpoints,
and the density update conserves mass.
Homogeneous Neumann conditions are used for the chemical variables.

For the correction, we assume $\varepsilon>0$, omit the superscript $n+1$,
and denote the provisional solution of
\eqref{eq:modified-main-scheme} by $\widehat\rho$ and $\widehat g$,
with $\langle\widehat g_{j,\cdot}\rangle_h=0$.
Define
\[
\widehat f_{j+\frac12,k}
=\widehat\rho_{j+\frac12}\psi_0(v_k)
+\frac{\varepsilon}{2}
 \bigl(\widehat g_{j,k}+\widehat g_{j+1,k}\bigr).
\]
If $\widehat\rho_{j+\frac12}\ge0$ and
$\widehat f_{j+\frac12,k}\ge0$ for all $j,k$,
we set $\rho=\widehat\rho$ and $g=\widehat g$
without further correction.
Otherwise, we apply the following procedure.
For $\widehat\rho_{j+\frac12}>0$, set
\[
\theta_{j+\frac12}
=\min\left\{1,\,
\min_{\widehat f_{j+\frac12,k}<0}
\frac{\widehat\rho_{j+\frac12}\psi_0(v_k)}
{\widehat\rho_{j+\frac12}\psi_0(v_k)-\widehat f_{j+\frac12,k}}
\right\},
\]
where the inner minimum is omitted if all values are nonnegative, and
\[
\rho_{j+\frac12}^{*}=\widehat\rho_{j+\frac12},
\qquad
\bar g_{j+\frac12,k}^{*}
=\frac{\theta_{j+\frac12}}{2}
 \bigl(\widehat g_{j,k}+\widehat g_{j+1,k}\bigr).
\]
For $\widehat\rho_{j+\frac12}\le0$, set
$\rho_{j+\frac12}^{*}=\bar g_{j+\frac12,k}^{*}=0$.
Then $\rho_{j+\frac12}^{*}\psi_0(v_k)
+\varepsilon\bar g_{j+\frac12,k}^{*}\ge0$ and
$\langle\bar g_{j+\frac12,\cdot}^{*}\rangle_h=0$.
Assuming $\Delta x\sum_j\widehat\rho_{j+\frac12}>0$, restore the
input mass by
\[
\alpha=\frac{\sum_j\widehat\rho_{j+\frac12}}
                 {\sum_j\rho_{j+\frac12}^{*}},
\qquad
\rho_{j+\frac12}=\alpha\rho_{j+\frac12}^{*},
\qquad
\bar g_{j+\frac12,k}=\alpha\bar g_{j+\frac12,k}^{*}.
\]

We return the microscopic variable to interfaces by
\[
\begin{aligned}
\widehat h_{j,k}
 &=\tfrac12\bigl(\bar g_{j-\frac12,k}+\bar g_{j+\frac12,k}\bigr),
 &&1\le j\le J-1,\\
\widehat h_{0,k}
 &=\tfrac12\bigl(\bar g_{\frac12,k}+\bar g_{\frac12,-k}\bigr),&
\widehat h_{J,k}
 &=\tfrac12\bigl(\bar g_{J-\frac12,k}+\bar g_{J-\frac12,-k}\bigr),
\end{aligned}
\]
and set $h_j=(I-\Pi_h)\widehat h_j$.
Since this reconstruction need not preserve positivity, we check
\[
F_{j+\frac12,k}
=\rho_{j+\frac12}\psi_0(v_k)
+\frac{\varepsilon}{2}(h_{j,k}+h_{j+1,k}).
\]
If all values are nonnegative, we set $g_j=h_j$.
Otherwise, define
\[
\begin{gathered}
q_j=\frac{\min(\rho_{j-\frac12},\rho_{j+\frac12})}{|V|_h}
\quad(1\le j\le J-1),\qquad
q_0=\frac{\rho_{\frac12}}{|V|_h},\qquad
q_J=\frac{\rho_{J-\frac12}}{|V|_h},\\
m_j=\min_k h_{j,k},\qquad
\eta_j=
\begin{cases}
\min\{1,q_j/(-\varepsilon m_j)\},&m_j<0,\\
1,&m_j\ge0,
\end{cases}
\qquad g_{j,k}=\eta_j h_{j,k}.
\end{gathered}
\]
In this case, $\varepsilon g_{j,k}\ge-q_j$, and hence
\[
f_{j+\frac12,k}
=\rho_{j+\frac12}\psi_0(v_k)
+\frac{\varepsilon}{2}(g_{j,k}+g_{j+1,k})
\ge\frac{\rho_{j+\frac12}}{|V|_h}
-\frac{q_j+q_{j+1}}{2}\ge0.
\]
The interface scaling preserves the zero velocity mean and endpoint
reflection and leaves the density unchanged.
Together with the conservative density update, this gives
\[
\Delta x\sum_j\rho_{j+\frac12}
=\Delta x\sum_j\widehat\rho_{j+\frac12}
=\Delta x\sum_j\rho_{j+\frac12}^{\,n}.
\]
At fixed grids, this additional interface scaling is inactive for
sufficiently small $\varepsilon$ if $\rho$ is bounded below by a positive
constant independent of $\varepsilon$ and $h$ remains bounded.

\bibliographystyle{abbrv}
\bibliography{bib.bib}

\end{document}